\documentclass{article}
\usepackage{latexsym,amsfonts,palatino,eulervm,color,amsmath,amssymb,amsthm,graphics,graphicx,comment,float,enumitem,authblk}
\usepackage[hidelinks]{hyperref}
\usepackage[hang,flushmargin]{footmisc}

\usepackage{xcolor}
\usepackage{soul}

\definecolor{myblue}{RGB}{0, 102, 204}  
\definecolor{myred}{HTML}{FF5733}       

\newtheorem{theorem}{Theorem}[section]
\newtheorem{lemma}[theorem]{Lemma}
\newtheorem{definition}[theorem]{Definition}

\def\F{\mathbb{F}}

\def\N{\mathbb{N}}

\def\rank{\mathrm{rank}}
\def\rowspace{\mathrm{rowspace}}

\providecommand{\keywords}[1]{\textbf{Keywords:} #1}
\providecommand{\MSC}[1]{\textbf{MSC:} #1}
 
\title{The rank distribution of matrices representing graphs with a long induced path over the field of two elements}
\author[1,2]{Badriah Safarji}
\author[3]{Cian O'Brien}
\author[2]{Rachel Quinlan}
\affil[1]{Department of Mathematics, Faculty of Sciences of Tabuk, University of Tabuk, 47512 Tabuk, Saudi Arabia}
\affil[2]{School of Mathematical and Statistical Sciences, University of Galway}
\affil[3]{Department of Mathematics and Computer Studies, Mary Immaculate College}

\begin{document}
\setlength{\parindent}{0pt}
\setlength{\parskip}{1ex}

\maketitle

\begin{abstract}
 
A square matrix $M$ represents a graph $\Gamma$ if its nonzero off-diagonal entries encode the adjacencies of $\Gamma$, subject to a fixed ordering of the vertices. Over the field of two elements, we investigate the distribution of ranks in the affine space consisting of all matrices representing a given $\Gamma$. In particular, we consider which graphs of order $n$ are represented by more matrices of rank $n-1$ than of rank $n$. This property reflects an exceptional feature of the space $M_n(\F_2)$ of all $n\times n$ matrices over $\F_2$, namely that its most frequently occurring rank is not $n$ but $n-1$. Our analysis focuses on the class of connected graphs with an induced path on all but one vertex. The main result is a characterisation of all such graphs that are represented by more matrices of rank $n-1$ than of rank $n$ over $\mathbb{F}_2$.

\end{abstract}


\MSC{05C50,  
15A03, 
15B33 
}

\section{Introduction}
Given a field $\mathbb{F}$ and a simple undirected graph $\Gamma$ with vertices $x_1,x_2,\dots,x_n$, a symmetric matrix $M$ with entries in $\mathbb{F}$ \emph{represents} $\Gamma$ (with respect to this ordering of the vertices) if, when $i\ne j$, the $(i,j)$-entry of $M$ is $0$ if and only if there is no edge between $x_i$ and $x_j$ in $\Gamma$. The diagonal entries of $M$ are not subject to any constraints, and therefore there are many matrices representing $\Gamma$ over $\F$. Since the non-zero off-diagonal entries may vary through the non-zero elements of $\F$, subject to the constraint of symmetry, and the diagonal entries are freely chosen from $\F$, the set $S(\F,\Gamma)$ of all symmetric matrices over $\F$ representing $\Gamma$ is a union of affine spaces of dimension $n$. If $r$ is the rank of some matrix $A\in S(\F,\Gamma )$ and $r<n$, then one may alter a single diagonal entry in $A$ to obtain an element of $S(\F,\Gamma )$ of rank $r+1$. It follows that the set of all ranks of matrices in $S (\Gamma)$ is an interval of the form $\{r\in\N :m\le r\le n\}$,where the positive integer $m$ is called the \emph{minimum rank} of $\Gamma$ over $\F$. The problem of determining $m$ is known as the \emph{minimum rank problem for graphs}. An informative survey of the extensive literature on the minimum rank problem (until 2007) is provided by Fallat and Hogben in \cite{fallat2007minimum}. The minimum rank generally depends on the choice of field $\F$. The case $\F=\mathbb{R}$ has been the subject of particularly concentrated and sustained attention, but the case of finite fields offers scope to determine more detailed information, due to the availability of combinatorial and enumerative methods. For every positive integer $k$ and prime power $q$, a description of the structure of all graphs whose minimum rank over $\F_q$ is at most $k$ is given in \cite{grout2008minimum}. The average minimum rank of all labelled graphs of order $n$ over a finite field is investigated in \cite{Friedland_Loewy}.  The relationship between the minimum rank of a simple graph over $\F_2$ and the minimum cardinality of a \emph{subgraph complementation system} of the graph is studied in \cite{MR4395941}. Graphs with minimum rank at most 3 over $\F_2$ are characterized in \cite{MR2489366}, via a list of 62 minimal forbidden subgraphs of order at most 8.

For the finite field $\F_q$ of order $q$, one can investigate the numbers of matrices of each rank in the finite set $S(\F_q,\Gamma )$. For $k\le n$, we write $R_{k,q}(\Gamma )$ for the number of matrices of rank $k$ in $S(\F_q,\Gamma )$.

\begin{definition}
   Let $\Gamma$ be a graph of order $n$. The \emph{rank distribution} of $\Gamma$ over the finite field $\F_q$ is the list $R_{m,q}(\Gamma),\dots ,R_{n,q}(\Gamma )$, where $m$ is the minimum rank of $\Gamma$ over $\F_q$. 
\end{definition}

 In this paper, we focus on rank distributions over $\F_2$. This is because $\F_2$ is the only finite field over which the list of numbers of $n \times n$ matrices of rank $0,1,2,\dots,n$ is not a strictly increasing sequence, as discussed below. The same is true when restricting to symmetric matrices.

Another reason that the field of two elements is exceptional in this context is that, subject to a fixed ordering of the vertices, the off-diagonal entries of the matrix representing a graph are fully determined. This means that the set of all such matrices is an affine subspace of $M_n(\F_2)$ of dimension $n$.

Fisher and Alexander \cite{fisher1966matrices} provide the following formula for the number $M(n,k,q)$ of $n\times n$ matrices of rank $k$ over the field $\mathbb{F}_q$, where $0 \leq k \leq n$.
\[M(n,k,q) = \frac{\prod\limits_{j=0}^{k-1}\left(q^n-q^j\right)^2}{\prod\limits_{j=0}^{k-1}\left(q^{k}-q^j\right)}\]

Comparing $M(n,k,q)$ and $M(n,k+1,q)$, we have the following.
\begin{equation*}
\begin{aligned}
   M(n,k+1,q) 
   &= \frac{\prod\limits_{j=0}^{k}\left(q^n-q^j\right)^2}{\prod\limits_{j=0}^{k}\left(q^{k+1}-q^j\right)} 
   = \frac{\left(q^n-q^k\right)^2\prod\limits_{j=0}^{k-1}\left(q^n-q^j\right)^2}{q^{k+1}\prod\limits_{j=0}^{k}\left(q^{k}-q^{j-1}\right)}
   = \frac{\left(q^n-q^k\right)^2}{q^{{k+1}}}\cdot \frac{\prod\limits_{j=0}^{k-1}\left(q^n-q^j\right)^2}{\prod\limits_{j=-1}^{k-1}\left(q^{k}-q^j\right)}\\
   &=\frac{q^{2k}\left(q^{n-k}-1\right)^2}{q^{k+1}\left(q^k-q^{-1}\right)}\cdot \frac{\prod\limits_{j=0}^{k-1}\left(q^n-q^j\right)^2}{\prod\limits_{j=0}^{k-1}\left(q^{k}-q^j\right)} = \frac{q^k\left(q^{n-k}-1\right)^2}{q^{k+1}-1}M(n,k,q)
\end{aligned}
\end{equation*}

$M(n,k,q)>M(n,k+1,q)$ if and only if 
$\frac{q^k\left(q^{n-k}-1\right)^2}{q^{k+1}-1} < 1$, so 
$q^{2n-k} - 2q^n -q^{k+1} + q^k + 1 < 0$.

If $k \le n-2$:
\begin{equation*}
\begin{aligned}
   q^{2n-k} - 2q^n -q^{k+1} + q^k + 1 &\ge q^{n+2} - 2q^n -q^{n-1} + q^{n-2} + 1\\
   &\ge q^{n-2}\left(q^4 - 2q^2-q+1\right) + 1 > 0
\end{aligned}
\end{equation*}

If $k = n-1$:
\begin{equation*}
\begin{aligned}
   q^{2n-k} - 2q^n -q^{k+1} + q^k + 1 &= q^{n+1} - 2q^n -q^{n} + q^{n-1} + 1\\
   & = q^{n-1}(q^2-3q+1)+1
\end{aligned}
\end{equation*}

This expression is positive for $q\ge 3$. However, if $q=2$, it is equal to $1-2^{n-1}$. For fixed $q\ge 3$ and fixed $n\ge 1$, we conclude that $M(n,k,q)$ strictly increases with $k$, for $0\le k\le n$. For fixed $n\ge 2$, $M(n,k,2)$ increases with $k$ for $0\le k\le n-1$ but $M(n,n-1,2)>M(n,n,2)$. Therefore, the most frequently occurring rank in $M_n(\F_2)$ is not $n$ but $n-1$. 


When we restrict attention to the symmetric matrices, $\F_2$ is again an exception. MacWilliams \cite[Theorem 2]{macwilliams1969orthogonal} provides the following formulae for the number $N(n, r, q)$ of symmetric $n \times n$ matrices over $\F_q$ of rank $r=2s$ and rank $r=2s+1$, $s \in \mathbb{N}$.
\begin{equation}\label{symmetric count even}
   N(n, 2s, q)=\prod_{i=1}^s \frac{q^{2 i}}{q^{2 i}-1} \cdot \prod_{i=0}^{2 s-1}\left(q^{n-i}-1\right), 2s \leq n
\end{equation}

\begin{equation}\label{symmetric count odd}
    N(n, 2s+1, q)=\prod_{i=1}^s \frac{q^{2i}}{q^{2i}-1} \cdot \prod_{i=0}^{2 s}\left(q^{n-i}-1\right),2 s+1 \leq n.
\end{equation}


Again we compare these counts for consecutive ranks.
 
\underline{Suppose $r = 2s+1 \le n$:}
\begin{equation*}
\begin{aligned}
   N(n,2s+1,q) &= \prod_{i=1}^s \frac{q^{2i}}{q^{2i}-1} \cdot \prod_{i=0}^{2 s}\left(q^{n-i}-1\right)\\
   &=\left(q^{n-2s}-1\right)\cdot  \prod_{i=1}^{s} \frac{q^{2 i}}{q^{2 i}-1} \cdot \prod_{i=0}^{2 s-1}\left(q^{n-i}-1\right)=\left(q^{n-2s}-1\right)\cdot  N(n,2s,q)\\
\end{aligned}
\end{equation*}

\underline{Suppose $r = 2s \le n$:}
\begin{equation*}
\begin{aligned}
   N(n,2s,q) &= \prod_{i=1}^s \frac{q^{2 i}}{q^{2 i}-1} \cdot \prod_{i=0}^{2 s-1}\left(q^{n-i}-1\right)\\
   &=\frac{q^{2s}}{q^{2s}-1} \cdot \left(q^{n-(2s-1)}-1\right)\cdot  \prod_{i=1}^{s-1} \frac{q^{2 i}}{q^{2 i}-1} \cdot \prod_{i=0}^{2 s-2}\left(q^{n-i}-1\right)\\
   &= \left(1+\frac{1}{q^{2s}-1} \right) \cdot \left(q^{n+1 - 2s} -1\right)\cdot  N(n,2s-1,q)\\
\end{aligned}
\end{equation*}

Setting $q=2$, we observe that $$N(n,n,2) = N(n,n-1,2)$$ for odd $n$, and $$N(n,n,2) = \left(1+\frac{1}{2^n-1}\right)N(n,n-1,2)$$ for even $n$. Otherwise, if $r < n$ or $q > 2$, then $N(n,r,q)$ is an integer ($\ge2$) multiple of $N(n,r-1,q)$. For this reason, we generally expect a connected graph of order $n$ to be represented by at least as many matrices over $\F_q$ of rank $r$ as of rank $r-1$. Exceptions to this pattern must be prevalent in the case $q=2$ and $r=n$.

The goal of this article is to identify classes of connected graphs which are represented by more matrices over $\F_2$ of rank $n-1$ than of rank $n$. We focus on connected graphs that have a path on all but one of their vertices as an induced subgraph. This is motivated by the special role of the path graph in the study of the minimum rank problem. Over any field $\F$, the path $P_n$ is the unique graph of order $n$ whose minimum rank over $\F$ is $n-1$ \cite{fallat2007minimum}.

To demonstrate the variety in rank distribution possible across different classes of graphs, we consider the examples of complete graphs, cycle graphs, and path graphs.  If $M$ is a matrix representing $K_n$ over $\F_2$ with $r$ zero entries on the diagonal, then the rank of $M$ is $r+1$ if $r<n$. If $r=n$, then $M$ has rank $n$ if $n$ is even, and rank $n-1$ if $n$ is odd.

So for $r \leq n-2$, the number of matrices of rank $r$ representing $K_n$ is the number of ways to put $r-1$ zeros on the diagonal, ${n \choose r-1}$. If $n$ is even, the number of rank $n-1$ is ${n \choose n-2}$ and of rank $n$ is ${n \choose n-1}+1$. If $n$ is odd, the number of rank $n-1$ is ${n \choose n-2}+1$ and of rank $n$ is ${n \choose n-1}$.

Theorem \ref{path numbers} describes the rank distribution of the path graph $P_n$, and Theorem \ref{cycle} describes the rank distribution of the cycle graph $C_n$. The table below summarises these results, along with the rank distribution of the complete graph $K_n$. We use $R_k(\Gamma)$ to denote the number of matrices over $\F_2$ of rank $k$ representing $\Gamma$.
\begin{center}
\renewcommand{\arraystretch}{1.5}
\begin{tabular}{|c||c|c|c|c|c|c|c|}
\hline
$\Gamma$ & $R_1(\Gamma)$ & $R_2(\Gamma)$ & $\dots$ & $R_{n-3}(\Gamma)$ & $R_{n-2}(\Gamma)$ & $R_{n-1}(\Gamma)$ & $R_{n}(\Gamma)$\\
\hline
\hline
$K_n$ & ${n \choose 0}$ & ${n \choose 1}$ & $\dots$ & ${n \choose n-4}$ & ${n \choose n-3}$ & ${n \choose n-2} + \frac{1-(-1)^n}{2}$ & ${n \choose n-1} + \frac{1+(-1)^n}{2}$\\
\hline
$C_n$ &  &  & $\dots$ & & $\frac{1}{3}\big(2^{n-1}+(-1)^n\big)$ & $2^{n-1}$ & $\frac{1}{3}\big(2^n+(-1)^{n-1}\big)$\\
\hline
$P_n$ &  &  & $\dots$ & &  & $\frac{1}{3}\big(2^n+(-1)^{n-1}\big)$ & $\frac{1}{3}\big(2^{n+1}+(-1)^n\big)$\\
\hline
\end{tabular}
\end{center}

So the most frequently occurring rank(s) in the rank distribution of $K_n$ is $\frac{n+2}{2}$ if $n$ is even, and $\frac{n+1}{2}$ and $\frac{n+3}{2}$ if $n$ is odd, while the most frequently occurring rank for $C_n$ is $n-1$, and for $P_n$ it is $n$.

In this work, we restrict our attention to connected graphs because the rank distribution of any graph can be determined by the rank distributions of its connected components. The vertices of a disconnected graph can be ordered such that each matrix representing it is block-diagonal with blocks corresponding to its connected components, and so the rank of a matrix representing a disconnected graph is the sum of the ranks of  the blocks corresponding to its components. Therefore, studying connected graphs captures all the essential complexity of the general case. 

We have seen that the number of symmetric $n \times n$ matrices over $\F_q$ of rank $r$ is at least double the number of rank $r-1$ for all $r \in \{1,\dots,n\}$ with $(q,r)\ne(2,n)$, while the number of rank $n-1$ and of rank $n$ are approximately equal over $\F_2$. Since every symmetric matrix represents some graph, there exist graphs of order $n$ represented by more matrices over $\F_2$ of rank $n-1$ than of rank $n$. Disconnected graphs alone do not account for this. There are more connected graphs of order $n$ than disconnected graphs because the complement of a disconnected graph is connected, but the complement of a connected graph is not necessarily disconnected. The rank distribution of a disconnected graph is equal to the discrete convolution of the distributions of its connected components. Consequently, the more connected components a graph has, the more likely its most frequently occurring rank is to be strictly between its minimum and maximum rank. This is true over any finite field, and therefore does not account for the exceptional rank distribution of all symmetric matrices over $\F_2$. 



Throughout this paper, we write $\mathcal{G}_d$ for the class of graphs containing an induced path on all but one of its vertices, and for which the vertex that does not belong to the induced path has degree $d \in \mathbb{N}$. The rest of the paper is arranged as follows. In Section \ref{path section}, we describe the rank distribution for the path graph $P_n$. In Section \ref{induced path section}, we define functions to count the number of matrices representing graphs of order $n$ in $\mathcal{G}_d$ over $\F_2$ of rank $n-1$ and $n$, and describe formulae for these functions in terms of nullspace vectors of matrices representing $P_n$. In Section \ref{reduction}, we define recurrence relations for these counting functions, expressing their value for graphs in $\mathcal{G}_d$ in terms of graphs in $\mathcal{G}_{d-1}$. In Section \ref{alpha}, we begin by classifying all graphs in $\mathcal{G}_1$ represented by more matrices of rank $n-1$ than $n$. We then use the results from $\mathcal{G}_1$ to determine all such graphs in $\mathcal{G}_2$, continuing similarly for $\mathcal{G}_3$, $\mathcal{G}_4$, and $\mathcal{G}_5$. Finally, we show that there are no such graphs in $\mathcal{G}_d$ for $d \ge 6$. These results are then summarised in Theorem \ref{recap}.

\section{The $\F_2$-rank distribution of the path graph}
\label{path section}
For a positive integer $n$, we write $P_n$ for the path graph on $n$ vertices. We write $x_1, x_2, \ldots, x_n$ for the vertices of $P_n$, where $x_1$ and $x_n$ are the two vertices of degree 1, and $x_i$ is adjacent to $x_{i+1}$ for $1 \leqslant i \leqslant n-1$. 
Over any field $\mathbb{F}$, every matrix that represents $P_n$ with respect to this vertex ordering is tridiagonal with nonzero entries in the super-diagonal and sub-diagonal; therefore, it has rank at least $n-1$. Over any field $\mathbb{F}$, it is routine to check that if the nonzero off-diagonal entries of a matrix representing $P_n$ are all 1, then the diagonal entries may be completed to obtain a matrix either of rank $n$ or of rank $n-1$. The path $P_n$ has a special role in the minimum rank problem for graphs; for $n \geqslant 1$ and for every field $\mathbb{F}, P_n$ is the unique graph whose minimum rank over $\mathbb{F}$ is $n-1$ (see \cite{fallat2007minimum}).

We define an \emph{indeterminate matrix} over $\mathbb{F}$ to be a matrix in which each entry is either an element of $\mathbb{F}$ or an indeterminate. A \emph{completion} of an indeterminate matrix $M$ over $\mathbb{F}$ is the matrix that results from an assignment of elements of $\mathbb{F}$ to the indeterminates of $M$. We say an indeterminate matrix $M$ \emph{represents} a graph $\Gamma$ if, for all $i\ne j$, $M_{ij} = 0$ when $x_ix_j$ is not an edge of $\Gamma$ and $M_{ij}=M_{ji}$ if $M_{ij}$ is an indeterminate. Therefore any completion of $M$ in which off-diagonal indeterminates are assigned non-zero values represents $\Gamma$.

We define $M(P_n)$ to be the indeterminate matrix over $\mathbb{F}_2$ with all entries in the first superdiagonal and the first subdiagonal equal to 1, indeterminates on the main diagonal, and zeros elsewhere.
$$
M(P_n) =
\begin{bmatrix}
d_1 & 1 & 0 & \cdots & 0 \\
1 & d_2 & 1 & \ddots & \vdots \\
0 & 1 & d_3 & \ddots & 0 \\
\vdots & \ddots & \ddots & \ddots & 1 \\
0 & \cdots & 0 & 1 & d_n
\end{bmatrix}
$$

Our goal in this section is to determine $R_n(P_n)$ and $R_{n-1}(P_n)$, the number of $\mathbb{F}_2$-completions of $M\left(P_n\right)$ of ranks $n$ and $n-1$, respectively. We denote the $i$th standard basis vector by $e_i$.

\begin{lemma}\label{path rank}
 Let $A$ be a completion of $M(P_n)$. Then $A$ has rank $n$ if and only if $e_1\in\F_2^{n}$ is in the column space of $A$. Equivalently, $A$ has rank $n$ if and only if $e_n\in\F_2^n$ is in the column space of $A$. 
\end{lemma}

\begin{proof} 
If $A$ has rank $n$, then its column space is $\mathbb{F}_2^n$ which includes $e_1$. On the other hand, suppose that $e_1$ belongs to the column space of $A$. Since the first column of $A$ has the form $a e_1+e_2$ for $a\in\F_2$, it follows that $e_2$ is also in the column space of $A$. Applying this reasoning to each successive column of $A$, we observe that all the standard basis vectors of $\mathbb{F}_2^n$ belong to the column space of $A$, so $A$ has rank $n$. The same argument applies to $e_n$, starting with column $n$.
\end{proof} 

The following theorem shows that approximately one-third of $\mathbb{F}_2$-matrices representing $P_n$ have rank $n-1$, with the remainder having rank $n$.

\begin{theorem}\label{path numbers}
$
R_n\left(P_n\right)=\displaystyle\frac{1}{3}\left(2^{n+1}+(-1)^n\right), R_{n-1}\left(P_n\right)=\displaystyle\frac{1}{3}\left(2^n+(-1)^{n+1}\right)
$
\end{theorem}
\begin{proof}
Let $A^{\prime}$ be a completion of $M(P_{n-1})$. Then $A'$ has rank $n-1$ or $n-2$. Let $A$ be the partial completion of $M(P_n)$ that has $A^{\prime}$ as its lower right $(n-1) \times(n-1)$ submatrix, and the indeterminate $d_1$ as its upper left entry.

First, suppose that $A^{\prime}$ has rank $n-2$. By Lemma \ref{path rank}, the first standard basis vector of $\mathbb{F}_2^{n-1}$ is not in the column space of $A'$, nor is its transpose in the row space of $A'$. It follows that the matrix consisting of the last $n-1$ columns of $A$ has rank $n-1$ and that both completions of $A$ to an element of $M_n(\mathbb{F}_2)$ have rank $n$, since the first column is independent of the remaining columns, regardless of the value assigned to $d_1$. So $A$ has rank $n$ for both choices of $d_1$, which means that every rank $n-2$ completion of $M(P_{n-1})$ corresponds to two rank $n$ completions of $M(P_n)$.

Now suppose that $A'$ has rank $n-1$. Deleting the first row of $A$ leaves an $(n-1)\times n$ matrix of rank $n-1$, whose right nullspace contains a unique non-zero vector $u \in \mathbb{F}_2^n$. The first entry of $u$ is 1, since the columns of $A^{\prime}$ are linearly independent. The first entry of $Au$ is $d_1+u_2$, and it
follows that one choice of a value of $d_1$ determines a completion of rank $n-1$ of $M(P_n)$ that has $u$ in its right nullspace, and the other determines a completion of rank $n$. Therefore, each rank $n-1$ completion of $M(P_{n-1})$ corresponds to one rank $n$ completion and one rank $n-1$ completion of $M(P_n)$. 

We conclude that
\begin{equation}
\label{R(P_n)}
   R_n(P_n)=2 \cdot R_{n-2}(P_{n-1}) +R_{n-1}(P_{n-1}), \hspace{1cm} R_{n-1}(P_n)=R_{n-1}(P_{n-1}).
\end{equation}
The result follows by induction on $n$, noting that $R_2\left(P_2\right)=3$ and $ R_1\left(P_2\right)=1$.


\end{proof}

For an integer $n\ge 0$, we define
$$
F(n)=\frac{1}{3}(2^{n+1}+(-1)^{n}).
$$

For $n\ge 1$, it follows from Theorem \ref{path numbers} that $F(n)$ is equal to the number of $\F_2$-matrices of rank $n$ that represent $P_n$, and the number of rank $n$ that represent $P_{n+1}$. We note the following properties of $F(n)$, which will be useful in Sections \ref{reduction} and \ref{alpha}. 

\begin{lemma}\label{last n-1 rows}\phantom{Empty space}
\begin{enumerate} 
\item $F(n)$ is the number of completions of the matrix consisting of the last $n-1$ rows of $M(P_{n+1})$, whose rowspace avoids $e_1^\top$.    
\item $F(n)$ is the number of completions of the matrix consisting of the first $n-1$ rows of $M(P_{n+1})$, whose rowspace avoids $e_n^\top$.    
\end{enumerate}
\end{lemma}

\begin{proof} 
Deleting the first row from a completion of $M(P_{n+1})$ of rank $n$ leaves an $n\times (n+1)$ matrix of rank $n$, whose rowspace avoids $e_1^\top$ by Lemma \ref{path rank}. On the other hand, let $A$ be a completion of the last $n$ rows of $M(P_{n+1})$ whose rowspace avoids $e_1^\top$. The rowspace of $A$ intersects $\langle e_1^\top,e_2^\top\rangle$ in a 1-dimensional subspace that contains exactly one of $e_2^\top$ and $e_1^\top +e_2^\top$. Hence, the insertion of an additional row at the top can extend $A$ in exactly one way to a completion of $M(P_{n+1})$ of rank $n$. These observations establish  a bijective correspondence that proves the first statement. The second is proved in a similar way.    
\end{proof}






\begin{lemma}  
\label{F(p)}  
For all positive integers $p$ and $q$, the function $F$ satisfies the following recurrence relations. 

\begin{enumerate}

\item $2 F(p-1)=F(p)+(-1)^{p+1}$.  
\item $4 F(q-1) F(p-1)=F(q) F(p)+(-1)^{p+1} F(q)+(-1)^{q+1} F(p)+(-1)^{p+q}$.

\end{enumerate}
\end{lemma}

\section{Graphs with a long induced path}\label{induced path section}

The remainder of this article is concerned with the class $\mathcal{G}$ of connected graphs that have a path on all but one of their vertices as an induced subgraph. The motivation for studying this class of graphs is provided by Theorem \ref{path numbers}, which fully describes the $\F_2$-rank distribution of $P_n$. For a graph $\Gamma$ of order $n$ in $\mathcal{G}$, we investigate the relationship between the rank of an $\F_2$-matrix representing $\Gamma$ and that of its submatrix corresponding to an induced subgraph isomorphic to $P_{n-1}$.

For any graph $\Gamma$ of order $n$ in $\mathcal{G}$, we list the vertices of $\Gamma$ as $x,x_1,\dots ,x_{n-1}$, where the subgraph induced on $\{x_1,\dots ,x_{n-1}\}$ is a path with edges $x_ix_{i-1}$ for $1\le i\le n-2$. 
We write $M(\Gamma )$ for the indeterminate matrix that generically represents $\Gamma$ with respect to this vertex ordering. Then 
$$
M(\Gamma)=\left(\begin{array}{c|c}
d_0 & v ^\top \\
\hline
v & M\left(P_{n-1}\right)
\end{array}\right),
$$
where the upper left entry $d_0$ is an indeterminate and the vector $v \in \mathbb{F}_2^{n-1}$ records the incidences at the vertex $x$. 
Since every $\F_2$-completion of $M\left(P_{n-1}\right)$ has rank at least $n-2$, every $\F_2$-matrix representing $\Gamma$ has one of three possible ranks: $n-2, n-1$ or $n$.

The following theorem details how the rank of a completion of $M(P_{n-1})$ determines the ranks of its two extensions to completions of $M(\Gamma )$.

\begin{theorem}
\label{rank n-2}
    Let $A$ be a completion of $M(P_{n-1})$ and let $A(0)$ and $A(1)$ be the completions of $M(\Gamma)$ respectively given by 
    $$
    A(0)=\left(\begin{array}{c|c}
0 & v^\top \\
\hline
v & A
\end{array}\right), \ \ \ \ A(1)=\left(\begin{array}{c|c}
1 & v ^\top \\
\hline
v & A
\end{array}\right).
    $$
    Then 
    \begin{enumerate}
        \item If $\rank (A)=n-1$, then one of $A(0)$ and $A(1)$ has rank $n-1$ and the other has rank $n$.
        \item If $\rank (A)=n-2$ and $v ^\top \not\in\rowspace (A)$, then both $A(0)$ and $A(1)$ have rank $n$.
        \item If $\rank (A)=n-2$ and $v ^\top\in\rowspace (A)$, then one of $A(0)$ and $A(1)$ has rank $n-2$ and the other has rank $n-1$.
    \end{enumerate}
\end{theorem}

\begin{proof}\renewcommand{\qedsymbol}{}
Let $M$ denote the indeterminate matrix obtained from $M(\Gamma )$ by completing $M(P_{n-1})$ to $A$, and retaining the indeterminate $d_0$ in the $(1,1)$ position.
\begin{enumerate}

\item
Suppose that $\rank (A)=n-1$. Then $A(0)$ and $A(1)$ both have rank at least $n-1$. Let $A'$ denote the $(n-1)\times n$ submatrix of $M$ consisting of rows 2 through $n$, which are linearly independent in $\F_2^n$. The rows of $A$ form a basis of $\F_2^{n-1}$, and there is a unique $w\in\F_2^{n-1}$ for which $w ^\top A =v ^\top$. Then $w^\top A'$ is either equal to the first row of $A(0)$ or of $A(1)$, and exactly one of $A(0)$ and $A(1)$ has rank $n-1$. The other has rank $n$, since its first row is not a linear combination of subsequent rows.

\item 
Now suppose that $\rank (A)=n-2$ and that $v ^\top$ is not in the rowspace of $A$. Then (since $A$ is symmetric) $v$ is not a linear combination of the columns of $A$. It follows that extending $A$ to either $A(0)$ or $A(1)$ increases the rank from $n-2$ to $n$.

\item 
Since $v ^\top$ is a linear combination of the rows of $A$, either $(0|v^\top)$ or $(1|v^\top)$ is a linear combination of the rows of the $(n-1)\times n$ matrix $(v|A)$, which has rank $n-2$. Hence at least one of $A(0)$ and $A(1)$ has rank $n-2$. Both have rank $n-2$ if and only if the transpose of $e_1\in\F_2^{n}$ belongs to the rowspace of $(v^\top  | A)$, which means $u ^\top(v|A)=e_1^\top$ for some $u\in\F_2^n$. This is impossible, since if $u^\top A=0$ then $u^\top v=0$ also, as $v$ is in the column space of $A$. \hfill $\square$ \end{enumerate}    \end{proof}

When the matrix $A$ of Theorem \ref{rank n-2} has rank $n-2$, there exists a unique nonzero right vector $u$ in the right nullspace of $A$. In this case, it is useful to determine whether the vector $v$ is orthogonal to $u$ or not.  If $u^\top v=0$, then $M(\Gamma)$ has one completion each of rank $n-1$ and $n-2$. If $u^\top v=1$, then both choices for $d_0$ result in a rank $n$ completion of $M(\Gamma)$. As a result, we restrict our attention to completions of $M(\Gamma)$ for which the lower-right $(n-1)\times(n-1)$ submatrix has rank $n-2$.

For a graph $\Gamma$ in the class $\mathcal{G}$, we write $A(\Gamma )$ and $B(\Gamma )$ respectively for the numbers of matrices of rank $n$ and $n-1$ that represent $\Gamma$ over $\F_2$ with respect to the vertex ordering 
$\{x,x_1,\dots ,x_{n-1}\}$, and for which the lower right $(n-1)\times (n-1)$ submatrix corresponding to the path on $x_1,\dots ,x_{n-1}$ has rank $n-2$. From Theorem \ref{rank n-2} if follows that 
$$
R_{n}(\Gamma )-R_{n-1}(\Gamma ) = A(\Gamma )-B(\Gamma ). 
$$
We write $\alpha (\Gamma) = A(\Gamma )-B(\Gamma )$, and proceed to identify those $\Gamma\in\mathcal{G}$ for which $\alpha (\Gamma)$ is negative.

\subsection{Vectors in the nullspace of completions of $M(P_n)$}

We now consider which column vectors over $\mathbb{F}_2$ may occur as the unique nonzero element of the right nullspace of a matrix representing $P_n$ in $M_n(\mathbb{F}_2)$.

\begin{lemma}
\label{nullspace}
		Suppose that $Au=0$, for a matrix $A\in M_n(\F_2)$ that represents $P_n$, and a nonzero column vector $u$. Then the first and last entries of $u$ are both 1, and $u$ has no pair of consecutive zero entries. 
\end{lemma}

\begin{proof}
We write $d_1,\dots ,d_n$ for the diagonal entries of $A$, and $u_1,\dots ,u_n$ for the entries of $u$. Since $Au=0,$ we have the following:
  
\begin{itemize}
\item $d_1 u_1+u_2=0$
\item $u_{i-1}+d_i u_i+u_{i+1}=0$, for $2\le i\le n-1$ 
\item $u_{n-1}+d_nu_n=0$
\end{itemize}
If $u_1=0$, then from the first of the above equations, it follows that $u_2=0$. Applying the second equation to the successive triples $(u_{i-1}, u_i,u_{i+1})$ from $i=2$, it follows that $u_i=0$ for all $i$. A similar argument applies if $u_n=0$, working through the triples in the opposite order. Thus, the zero vector is the only vector with the first or last entry equal to zero in the nullspace of any matrix that represents $P_n$. 
		
Suppose now that $u_i=u_{i+1}=0$ for some $i$ with $2\le i\le n-2$. Then $u_{i-1}$ and $u_{i+2}$ are also equal to zero, since $u_{i-1}+d_i u_i+u_{i+1}=0$ and $u_{i}+d_{i+1} u_{i+1}+u_{i+2}=0$. Repeating this argument, it follows that $u=0$. 
	\end{proof}
	
For a positive integer $n$, we write $U_ n$ for the set of vectors in $\F_2^ n$ that have no consecutive zero entries and have first and last entries equal to 1. Lemma \ref{nullspace} shows that every non-zero vector that is in the nullspace of a completion of $M(P_n)$ belongs to $U_n$. In the next lemma, we show that $U_n$ is exactly the set of non-zero vectors that occur in the nullspace of some rank $n-1$ completion of $M(P_n)$, and determine the number of such completions with a particular 1-dimensional nullspace.

\begin{lemma}\label{2^k}
Let $u\in U_{n}$ and let $z(u)$ be the number of zero entries in $u$. Then the number of matrices $A$ that represent $P_n$ and satisfy $Au=0$ is $2^{z(u)}$. 
\end{lemma}
	
\begin{proof}
We write $u_1,\dots ,u_n$ for the entries of $u$, and note that $u_1=u_n=1$. Let $A$ be a completion of $M(P_n)$, and write $d_1,\dots,d_n$ for the diagonal entries of $A$. Then $Au=0$ if and only if the following conditions are satisfied. 
\begin{itemize} 
\item $d_1+u_2=0$
\item For $2\le i\le n-1$, $u_{i-1}+d_iu_i+u_{i+1}=0$
\item $u_{n-1}+d_{n}=0$
\end{itemize}
The first and last equations above are satisfied if (and only if) $d_1=u_2$ and $d_n=u_{n-1}$. For $2\le i\le n-1$, we consider the cases $u_i=1$ and $u_i=0$ separately. If $u_i=1$, then $u_{i-1}+d_i u_i+ u_{i+1}=0$ is satisfied only by $d_i=u_{i-1}+u_{i+1}$. If $u_i=0$, then $u_{i-1}+u_{i+1}=1+1=0$, and the equation is satisfied  by both $d_i=0$ and $d_i=1$. It follows that every element of $U_{n}$ belongs to the nullspace of \emph{some} completion of $M(P_n)$ of rank $n-1$, and that for a particular $u\in U_{ n}$ with $z(u)$ zero entries, the number of choices for the diagonal entries of such a completion is $2^{z(u)}$.
\end{proof}


Lemma \ref{2^k} allows us to characterise $A(\Gamma)$ and $B(\Gamma)$ in the following way. For column vectors in $\F_2^n$, we write $\perp$ for the relation of orthogonality with respect to the standard scalar product.

\begin{theorem}
\label{A(G)}
    Let $\Gamma \in \mathcal{G}$ have order $n$, and let $v$ be the vector in $\F_2^{n-1}$ consisting of the last $n-1$ entries of the first column of $M(\Gamma )$. 
    Let $z(u)$ be the number of zero entries in $u$. Then
    \begin{enumerate}
        \item 
         $A(\Gamma ) = \sum \limits_{\substack{ u \in U_{n-1}\\u\not\perp v        }} 2^{z(u)+1}$ 
        \item
        $B(\Gamma ) = \sum \limits_{ \substack{ u \in U_{n-1}\\ u\perp v}} 2^{z(u)}$
        \end{enumerate}
\end{theorem}

\begin{proof}
Let $u \in U_{n-1}$, and let $T(u)$ be the set of completions of $M(\Gamma)$ whose lower right $(n-1)\times (n-1)$ submatrix is a completion of $M(P_{n-1})$ with nullspace $\langle u\rangle$. By Lemma \ref{2^k}, $u$ is in the nullspace of $2^{z(u)}$ rank $n-2$ completions of $M(P_{n-1})$, each of which extends to two elements of $T(u)$ through the choice of a value for $d_0$. Hence $|T(u)|=2^{z(u)+1}$. Note that $u\perp v \iff v^\top u = 0$, meaning $u\perp v$ if and only if $v^\top$ is in the rowspace of every rank $n-2$  completion of $M(P_{n-1})$ whose nullspace contains $u$.

If $u \not\perp v$, then it follows from Theorem \ref{rank n-2} that every element of $T(u)$ has rank $n$. If $u \perp v$, then it follows from Theorem \ref{rank n-2} that $2^{z(u)}$ elements of $T(u)$ have rank $n-1$ and $2^{z(u)}$ have rank $n-2$.  Together, these imply the result.
\end{proof}

\subsection{The $\F_2$-rank distribution of the cycle graph}

For any positive integer $k$, we write $T_k$ for the set of completions of $M(P_k)$ of rank $k-1$.

Let $v\in\F_2^{n-1}$. Then $v$ determines a graph $\Gamma (v)$ of order $n$ in $\mathcal{G}$. The vertex set of $\Gamma (v)$ is $\{x,x_1,\dots ,x_{n-1}\}$. Its edge set consists of the edges of an induced path $x_1x_2\dots x_{n-1}$, and those edges $xx_i$ for which $v_i=1$. We write $M(v)$ for the indeterminate matrix $M(\Gamma(v))$, with respect to the above ordering of the vertices, and we write $A(v),B(v)$ and $\alpha(v)$ respectively for the quantities $A(\Gamma (v)), B(\Gamma(v))$ and $\alpha(\Gamma(v))$.

For $T\in T_{n-1}$, let $u_T$ denote the unique non-zero element of the right nullspace of $T$. Then $v$ belongs to the column space of $T$ if and only if $v^\top u_T=0$. We define $T_0(v) = \{T\in T_{n-1}: v^\top u_T=0\}$ and $T_1(v)=\{T\in T_{n-1}(v):v^\top u_T =1\}$, and note that $T_{n-1}$ is the disjoint union of $T_0(v)$ and $T_1(v)$. Thus $|T_0(v)|+|T_1(v)|=|T_{n-1}|=F(n)$. Furthermore
$$
A(v)=2|T_1(v)|,\ B(v)=|T_0(v)|,\ \alpha(v)=A(v)-B(v)=2|T_1(v)|-|T_0(v)|=R_n(\Gamma (v))-R_{n-1}(\Gamma (v)). 
$$
Let $u\in U_{n-1}$. As noted in Theorem \ref{A(G)}, $u$ contributes $2^{z(u)+1}$ to $A(v)$ if $v^\top u=1$, and $u$ contributes $2^{z(u)}$ to $B(v)$ if $v^\top u=0$. If we fix $u$ and allow $v$ to vary through all the elements of $\F_2^{n-1}$, then the cases $v^\top u=0$ and $v^\top u=1$ occur with the same frequency, and the overall contribution of $u$ to $\sum_{v\in\F_2^n}A(v)$ exceeds its contribution to $\sum_{v\in\F_2^n}B(v)$ by a factor of 2. In this sense, we expect that $A(v)$ exceeds $B(v)$ by a factor of 2 ``on average", and that $\alpha (v)$ is positive on average. On the other hand, the elements of $U_{n-1}$ have the special form described in Lemma \ref{nullspace}; they are exactly those vectors in $\F_2^{n-1}$ that have no consecutive zero entries and have first and last entries equal to 1. For a given $v$, we do not expect the cases $v^\top u=0$ and $v^\top u=1$ to occur with the same frequency as $u$ ranges through $U_{n-1}$. There are two extreme cases. If $v \in \{e_1,e_{n-1}\}$, then $v^\top u =1$ for every $u\in U_{n-1}$. In this case, $\Gamma (v)$ is the path on $n$ vertices, $B(v)=0$, and $\alpha(v)$ is maximal. If $v=e_1+e_{n-1}$, then $v^\top u=0$ for all $u\in U_{n-1}$. In this case, $\Gamma (v)$ is the cycle on $n$ vertices, $A(v)=0$, and $\alpha(v)=-B(v)$ is minimal.

The next theorem gives a complete description of the $\F_2$-rank distribution of the cycle $C_n$.

\begin{theorem}\label{cycle}
$R_{n}(C_n)=\frac{1}{3}(2^n+(-1)^{n+1}),\ R_{n-1}(C_n)=2^{n-1},\ R_{n-2}(C_n)=\frac{1}{3}(2^{n-1}+(-1)^n)$    
\end{theorem}

\begin{proof}
    Let $A'$ be a completion of $M(P_{n-1})$, and let $M'$ be a matrix obtained from $M(C_n)$ by completing the last $n-1$ indeterminate so that the lower right of the submatrix is $A'$. The vector $v$ consists of the last $n-1$ entries of the first row of $M'$, which has 1 as its first and last entries and otherwise consists of zeros. Every element of $U_{n-1}$ is orthogonal to $v$. 

    If $\rank(A')=n-1$, then $v$ is in the rowspace of $A'$. By Lemma \ref{rank n-2}, one choice for the upper left entry gives a completion of $M(C_n)$ of rank $n-1$, and the other gives a completion of rank $n$. 

    If $\rank(A')=n-2$, then $v$ is again in the rowspace of $A'$ since $v$ is orthogonal to the element of $U_{n-1}$ that spans the nullspace of $A'$. By Lemma \ref{rank n-2}, one choice for the upper left entry gives a completion of $M(C_n)$ of rank $n-1$, and the other gives a completion of rank $n-2$. 

    So every matrix that represents $P_{n-1}$ and has rank $n-1$ contributes one to both $R_n(C_n)$ and $R_{n-1}(C_n)$, and every matrix that represents $P_{n-1}$ and has rank $n-2$ contributes one to both $R_{n-1}(C_n)$ and $R_{n-2}(C_n)$.
    From Theorem \ref{path numbers}, we conclude
    \begin{itemize}
    \item $R_n(C_n)=R_n(P_{n-1})=\frac{1}{3}(2^n+(-1)^{n+1})$.
    \item $R_{n-1}(C_n)=R_n(P_{n-1})+R_{n-1}(P_{n-1})=2^{n-1}$.
    \item $R_{n-2}(C_n)=R_{n-1}(P_{n-1})=\frac{1}{3}(2^{n-1}+(-1)^n)$
\end{itemize}
\end{proof}
Thus half of all $\F_{2}$-matrices representing the cycle $C_n$ have rank $n-1$, approximately one-third have rank $n$, and approximately one-sixth have rank $n-2$.

\section{Recurrences for $A(\Gamma)$ and $B(\Gamma)$}\label{reduction}

Let $\Gamma\in\mathcal{G}$, with vertices $x,x_1,\dots ,x_{n-1}$. Let $p$ be the minimal non-negative integer with $x_{p+1}$ adjacent to $x$, so that $p$ is the number of edges in the path $x_1\dots x_{p+1}$. Let $\Gamma_1$ be the graph obtained from $\Gamma$ by deleting $x_1,\dots ,x_p$ and their incident edges, and deleting the edge $xx_{p+1}$. Let $\Gamma_2$ be the graph obtained from $\Gamma_1$ by deleting the vertex $x_{p+1}$ and its incident edge. In this section we establish expressions for $A(\Gamma )$ and $B(\Gamma )$ in terms of the corresponding quantities for $\Gamma_1$ and $\Gamma_2$.

The degree of the vertex $x$ in both $\Gamma_1$ and $\Gamma_2$ is $\deg_\Gamma (x)-1$, allowing for recursive analysis of $\alpha (\Gamma)$ in terms of the degree of the vertex $x$ for graphs in $\mathcal{G}$.

\begin{theorem}\label{ABreduction}
Let $\Gamma\in\mathcal{G}$. Then for $\Gamma_1$, $\Gamma_2$, and $p$ defined as above:
\begin{enumerate}
    \item $A(\Gamma )=2F(p)B(\Gamma_1)+2F(p-1)A(\Gamma_2).$
        \item $B(\Gamma )=\frac{1}{2}F(p)A(\Gamma_1)+2F(p-1)B(\Gamma_2).$
\end{enumerate}
\end{theorem}
The proof of Theorem \ref{ABreduction} is presented in a series of steps. We recall that $A(\Gamma)$ and $B(\Gamma)$ are respectively the numbers of completions of $M(\Gamma )$ of rank $n$ and rank $n-1$, in which the lower right $(n-1)\times (n-1)$ submatrix is a completion of $M(P_{n-1})$ of rank $n-2$. We write $v$ for the vector in $\F_2^{n-1}$ consisting of the last $n-1$ entries of the first column of $M(\Gamma )$, which records the neighbours of $x$ in $\Gamma$. We note that the first $p$ entries of $v$ are zeros, and the first nonzero entry of $v$ is in position $p+1$.

We write $C(v)$ for the set of completions of rank $n-2$ of $M(P_{n-1})$ whose columnspace contains $v$, and $\overline{C}(v)$ for the set of rank $n-2$ completions of $M(P_{n-1})$ whose columnspace excludes $v$. Suppose that $M'$ is a completion of $M(\Gamma )$ that contributes either to $A(\Gamma )$ or $B(\Gamma )$, and let $M$ be the corresponding completion of $M(P_{n-1})$, which has rank $n-2$. Theorem \ref{rank n-2} implies that $M$ contributes to $A(\Gamma)$ if $M\in\overline{C}(v)$ and to $B(\Gamma )$ if $M\in C(v)$. Every $M\in\overline{C}(v)$ extends in two ways to a rank $n$ completion of $M(\Gamma)$, since both choices for the upper left entry result in matrices of rank $n$. However, every $M\in C(v)$ extends in only one way to a rank $n-1$ completion of $M(\Gamma )$, since the two choices for the upper left entry result in one matrix of rank $n-1$ and one of rank $n-2$. Hence 
\begin{equation}\label{A}
A(\Gamma ) = 2|\overline{C}(v)|,\ \ B(\Gamma )=|C(v)|.
\end{equation}
To prove Theorem \ref{ABreduction}, we need to express $|C(v)|$ and $|\overline{C}(v)|$ in terms of $p$ and the graphs $\Gamma_1$ and $\Gamma_2$. Each element $M$ of $C(v)$ or $\overline{C} (v)$ has a unique nonzero vector $u_M$ in its right nullspace. The entry $u_M[p+1]$ in position $p+1$ of $u_M$ is either 1 or 0. We define
\begin{eqnarray*}
{C}_1(v)=\{M\in {C}(v): u_M[p+1]=1\},& & C_0(v)=\{M\in C(v): u_M[p+1]=0\} \\
\overline{C}_1(v)=\{M\in\overline{C}(v): u_M[p+1]=1\},& & \overline{C}_0(v)=\{M\in\overline{C}(v): u_M[p+1]=0\}.
\end{eqnarray*}
The proof of Theorem \ref{ABreduction} depends on enumerating the elements of the four sets above in terms of $\Gamma_1$ and $\Gamma_2$. 

For a positive integer $t\ge 2$, we define $S(t)$ to be the $(t-1)\times t$ indeterminate matrix obtained from $M(P_t)$ by deleting the first row. Every $\F_2$-completion of $S(t)$ has rank $t-1$ and has a 1-dimensional right nullspace in $\F_2^{t}$. For a vector $v\in\F_2$, we define $\Gamma (v)$ to be the graph on $t+1$ vertices whose indeterminate matrix is 
$$
\left[\begin{array}{cc}
* & v^\top  \\ v & M(P_t)
\end{array}
\right].
$$
We denote this matrix by $M(v)$. \\

\begin{lemma}\label{S(t)}
For $t\ge 2$ and $v\in\F_2^t$, 
\begin{enumerate} 
\item $B(\Gamma (v))$ is the number of completions of $S(t)$ whose rowspace includes $v^\top$ and not $e_1^\top$.
\item $\frac{1}{2}A(\Gamma (v))$ is the number of completions of $S(t)$ whose rowspace includes neither $v^\top$ nor $e_1^\top$.
\end{enumerate}
\end{lemma}

\begin{proof}
The matrix $M(v)$ has $M(P_t)$ as its lower right $t\times t$ submatrix and $S(t)$ as its lower right $(t-1)\times t$ submatrix. Let $S$ be a completion of $S(t)$. If the rowspace of $S$ includes $e_1$, then it follows from Lemma \ref{path rank} that both extensions of $S$ to completions of $M(P_t)$ have rank $t$. Since $A$ and $B$ count completions with lower-right $t \times t$ submatrices of rank $t-1$, the extensions of these matrices to completions of $M(v)$ do not contribute to either $A(\Gamma (v))$ or $B(\Gamma (v))$.

If the rowspace of $S$ does not include $e_1^\top$, then $S$ has a unique extension to a completion $S'$ of rank $t-1$ of $M(P_t)$, whose rowspace is equal to that of $S$. The two extensions of $S'$ to completions of $\Gamma (v)$, determined by assigning a value from $\F_2$ to the upper left entry, potentially contribute to $A(\Gamma  (v))$ or $B(\Gamma (v))$. If $v$ does not belong to the rowspace of $S'$ (or equivalently $S$), then Theorem \ref{rank n-2} implies that both of these extensions have rank $t+1$, and both contribute to $A(\Gamma (v))$. On the other hand, if $v$ belongs to the rowspace of $S'$, then Theorem \ref{rank n-2} implies that one extension of $S'$ to a completion of $S$ has rank $t-1$ and the other has rank $t$, the latter of which contributes to $B(\Gamma (v))$. 

By Lemma \ref{path rank}, every completion of $M(v)$ that is counted by either $A(\Gamma (v))$ or $B(\Gamma (v))$ has a lower right $(t-1)\times t$ submatrix whose rowspace excludes $e_1^\top$. Therefore among the completions of $S$ whose rowspace excludes $e_1^\top$, the number whose rowspace includes $v^\top$ is $B(\Gamma (v))$, and the number whose rowspace excludes $v^\top$ is $\frac{1}{2}A(\Gamma(v))$.
\end{proof}

\begin{lemma}\label{0 in position p+1}
$|C_0(v)|=2F(p-1)B(\Gamma_2)$ and $|\overline{C}_0(v)|=F(p-1)A(\Gamma_2)$.
\end{lemma}

\begin{proof}

Let $M$ be a completion of $M(P_{n-1})$, and let $L$ and $R$ respectively denote its upper left $p\times p$ submatrix and its lower right $(n-p-2)\times (n-p-2)$ submatrices, which are respectively completions of $M(P_p)$ and $M(P_{n-p-2})$. Suppose that $M\in C_0(v)\cup\overline{C}_0(v)$. The right nullspace of $M$ contains a unique non-zero element $u$, with $u[p+1]=0$ and $u[p]=u[p+2]=1$. Then $Lu_1=0$ and $Ru_2=0$, where $u_1$ and $u_2$ are respectively the elements of $\F_2^p$ and $\F_2^{n-p-2}$ consisting of the first $p$ and the last $n-p-2$ components of $u$. We define $v_1$ and $v_2$ similarly. Since neither $u_1$ nor $u_2$ is the zero vector, it follows that $L$ and $R$ are both rank deficient. The vector $v^\top$ belongs to the rowspace of $M$ if and only if $v_2^\top u_2=0$, which occurs if and only if $v_2^\top$ is in the rowspace of $R$.

On the other hand, let $L$ be any completion of rank $p-1$ of $M(P_p)$ and let $R$ be any completion of rank $n-p-3$ of $M(P_{n-p-2})$. Let $u_1$ and $u_2$ be the non-zero elements of the right nullspaces of $L$ and $R$ respectively, and note that the last entry of $u_1$ and the first entry of $u_2$ are both 1. If the upper left $p\times p$ submatrix of $M(P_{n-1})$ is completed to $L$ and the lower right $(n-p-2)\times (n-p-2)$ region is completed to $R$, then both assignments of a value to the indeterminate in row $p+1$ result in a matrix of rank $n-2$, whose right nullspace contains the vector $u=[u_1 \ 0 \ u_2]$. Hence every choice for $L$ and $R$ contributes twice either to $|C_0(v)|$ or to $|\overline{C}_0(v)|$, according to whether $v_2$ belongs to the rowspace of $R$ or not. The number of choices for $L$ is $R_{p-1}(P_p)=F(p-1)$. Since $\Gamma_2 = \Gamma (v_2)$, the number of choices for $R$ with $v_2$ in its rowspace is $B(\Gamma_2)$, and the number of choices for $R$ with rowspace excluding $v_2$ is $\frac{1}{2}A(\Gamma_2)$. Hence $|C_0(v)|=2\times F(p-1)B(\Gamma_2)$ and $|\overline{C}_0(v)| = 2\times F(p-1)\times \frac{1}{2}A(\Gamma_2)=F(p-1)A(\Gamma_2)$.
\end{proof}

\begin{lemma}\label{1 in position p+1}
$|C_1(v)|=\frac{1}{2}F(p)A(\Gamma_1)$ and $|\overline{C}_1 (v)|=F(p)B(\Gamma_1)$.
\end{lemma}

\begin{proof}
Let $M\in C_1(v)\cup\overline{C}_1(v)$, so $M$ is a completion of $M(P_{n-1})$ of rank $n-2$, whose right nullspace includes a single non-zero vector $u$ with 1 in position $p+1$.
Then $e_{p+1}^\top$ is not in the rowspace of $M$, since $e_{p+1}^\top u=1$. The vector $v^\top$ belongs to the rowspace of $M$ if and only if $v^\top u=0$, which occurs if and only if $v_2^\top u_2=0$, where $v_2$ and $u_2$ are the vectors in $\F_2^{n-1-p}$ respectively consisting of the last $n-1-p$ entries of $v$ and of $u$. Let $R$ be the submatrix of $M$ in rows $p+2$ through $n-1$ and columns $p+1$ through $n-1$. Then $R$ has rank $n-p-2$ and its rowspace comprises exactly those vectors $w^\top$ with $w^\top u_2=0$. It follows that $v^\top$ belongs to the rowspace of $M$ if and only if $v_2^\top$ belongs to the rowspace of $R$. Let $v_2'$ be the element $v_2+e_1$ of $\F_2^{n-p-1}$, which differs from $v_2$ only in its first entry, which is 0. Then $(v_2')^\top$ is in the rowspace of $R$ if and only if $v_2^\top$ is not. Thus $M\in C_1(v)$ if and only if $(v_2')^\top$ is not in the rowspace of $R$. Alternatively $M\in\overline{C}_1(v)$, which occurs if and only if $(v_2')^\top$ is in the rowspace of $R$. Since $\Gamma_1=\Gamma (v_2')$, if follows from Lemma \ref{S(t)} that the number of possibilities for $R$ in an element of $C_1(v)$ or $\overline{C}_1(v)$ are respectively bounded above by $B(\Gamma_1)$ and $\frac{1}{2}A(\Gamma_1)$. If $L$ is the upper left $p\times (p+1)$ submatrix of $M$, then the rowspace of $L$ excludes $e_{p+1}^\top$ and $L$ occurs as the first $p$ rows of a unique completion of $M(P_{p+1})$ of rank $p$. It follows that the number of possibilities for $L$ is at most $R_p(P_{p+1})=F(p)$.

On the other hand, let $L'$ be a completion of the first $p$ rows of $M(P_{p+1})$ whose rowspace does not contain $e_{p+1}$, and let $R'$ be a completion of $S(n-p-1)$ whose rowspace does not contain $e_1^\top$. Let the non-zero vectors in the right nullspaces of $L'$ and $R'$ be $u_1$ and $u_2$ respectively. Noting that the last entry of $u_1$ and the first entry of $u_2$ are both 1,let $u$ be the vector in $\F_2^{n-1}$ that coincides with $u_1$ in its first $p+1$ entries and with $u_2$ in its last $n-1-p$ entries. Completing the upper left region of $M(P_{n-1})$ to $L'$, and the lower right region to $R'$ leaves a single choice for the indeterminate in row $p+1$, to ensures the resulting matrix $M'$ satisfies $M'u=0$ and hence belongs to $C_1(v)\cup\overline{C}_1(v)$. Thus the pair $(L',R')$ determines a unique element $M'$ of $C_1(v)\cup\overline{C}_1(v)$, which belongs to $C_1(v)$ if $(v_2')^\top$ is in the rowspace of $R'$ and to $\overline{C}_1(v)$ otherwise. Since the number of possibilities for $L'$ is $F(p)$, the conclusion follows from Lemma \ref{S(t)}.
\end{proof}

We complete the proof of Theorem \ref{ABreduction} by noting that
\begin{eqnarray*}
A(\Gamma ) & = &  2|\overline{C}(v)|=2|\overline{C}_1(v)|+2|\overline{C}_0(v)| = 2F(p)B(\Gamma_1)+2F(p-1)A(\Gamma_2), \\
B(\Gamma ) & = & |C(v)|=|C_1(v)|+|C_0(v)|=\frac{1}{2}F(p)A(\Gamma_1)+2F(p-1)B(\Gamma_2).
\end{eqnarray*}

\section{Characterising all graphs $\Gamma \in\mathcal{G}$ for which $\alpha(\Gamma)<0$}\label{alpha}

Recall $\alpha(\Gamma) = A(\Gamma)-B(\Gamma)$ counts the difference in the number of matrices of rank $n$ and $n-1$ representing $\Gamma\in\mathcal{G}$ over $\F_2$. Using Theorem \ref{ABreduction}, we derive the following expression for \(\alpha(\Gamma)\).
\[
\begin{aligned}
\alpha(\Gamma) &= 2F(p) B(\Gamma_1) + 2F(p - 1) A(\Gamma_2)  - \frac{1}{2}F(p) A(\Gamma_1) - 2F(p - 1) B(\Gamma_2) \\
&= \frac{1}{2}F(p)\big(4B(\Gamma_1) - A(\Gamma_1)\big) + 2F(p - 1)\big(A(\Gamma_2) - B(\Gamma_2)\big)
\end{aligned}
\]

To simplify the recurrence, we replace $A(\Gamma_2) - B(\Gamma_2)$ with $\alpha(\Gamma_2)$ and define $\beta(\Gamma) = 4B(\Gamma) - A(\Gamma)$.
The recurrence relation then becomes $
\alpha(\Gamma) = \frac{1}{2}F(p)\beta(\Gamma_1) + 2F(p - 1)\alpha(\Gamma_2)$.

By Theorem \ref{ABreduction}, we derive the following recurrence formula for \(\beta(\Gamma)\).
\[
\begin{aligned}
\beta(\Gamma) & =4B(\Gamma)-A(\Gamma ) \\
&= 4\left(\frac{1}{2}F(p) A(\Gamma_1) + 2F(p - 1) B(\Gamma_2)\right) - \left(2F(p) B(\Gamma_1) + 2F(p - 1) A(\Gamma_2)\right) \\
&= 2F(p)\left(A(\Gamma_1) - B(\Gamma_1)\right) + 2F(p - 1)\left(4B(\Gamma_2) - A(\Gamma_2)\right) \\
&= 2F(p)\alpha(\Gamma_1) + 2F(p - 1)\beta(\Gamma_2)
\end{aligned}
\]

We summarise these recurrences in the following lemma.
\begin{lemma}\phantom{hello}
\label{alpha-beta}\begin{enumerate}
\item $\alpha(\Gamma) = \frac{1}{2}F(p)\beta(\Gamma_1) + 2F(p - 1)\alpha(\Gamma_2)$
\item $\beta(\Gamma) = 2F(p)\alpha(\Gamma_1) + 2F(p- 1)\beta(\Gamma_2)$
\end{enumerate}
\end{lemma}

In this section, we use these recurrences to determine all $\Gamma \in \mathcal{G}$ with $\alpha(\Gamma) < 0$. For an integer $d\ge 1$, we write $\mathcal{G}_d$ for the class of graphs in $\mathcal{G}$ in which the vertex that does not belong to the induced path has degree $d$. If $\Gamma\in\mathcal{G}_d$ for $d>1$, the $\Gamma_1$ and $\Gamma_2$ belong to $\mathcal{G}_{d-1}$. The recurrence relations above express the values of $\alpha$ and $\beta$ for a graph in $\mathcal{G}_d$ in terms of corresponding values for graphs in $\mathcal{G}_{d-1}$.

\begin{figure}[h!]
    \centering
    \includegraphics[width=0.7\linewidth]{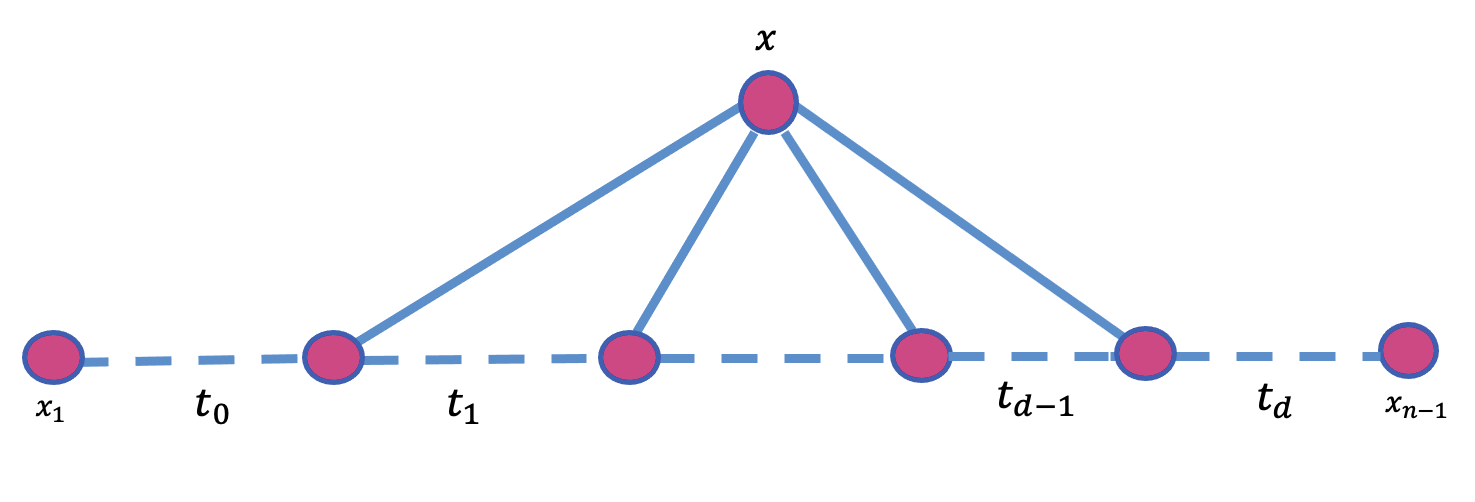}
    \vspace{-0.5cm}\caption{$\Gamma(t_0,t_1,\cdots,t_d)$}
    \label{Gamma(t)}
\end{figure}

 We write $\Gamma\left(t_{0}, t_1, \ldots, t_d\right)$ for the graph in $\mathcal{G}_d$ with the following properties, where $x$ is the vertex not in the induced path $P = x_1 x_2\dots x_{n-1}$, and the $i^\text{th}$ neighbour of $x$ is the neighbour of $i^\text{th}$ least index in $P$.
\begin{itemize}
    \item There are $t_0 \geq 0$ edges of $P$ between $x_1$ and the first neighbour of $x$.
    \item There are $t_i \geq 1$ edges of $P$ between the $i^\text{th}$ and $(i+1)^\text{th}$ neighbours of $x$ for $i \in \{1,\dots,d-1\}$.
    \item There are $t_d \geq 0$ edges of $P$ between the last neighbour of $x$ and $x_{n-1}$.
\end{itemize}

For $\Gamma = \Gamma(t_{0}, t_1, \ldots, t_d)$, we write $\alpha(t_{0}, t_1, \ldots, t_d) = \alpha(\Gamma)$, and similar for $\beta$, $A$, and $B$.

The function \( \beta \) plays an important role in determining whether \( \alpha(\Gamma) \) is negative for given $\Gamma\in\mathcal{G}$, since Lemma \ref{alpha-beta} implies that \( \alpha(\Gamma) \) can only be negative if either \( \beta(\Gamma_1) < 0 \) or \( \alpha(\Gamma_2) < 0 \). This allows us to determine which $\Gamma \in \mathcal{G}_d$ have $\alpha(\Gamma)<0$ solely in terms of graphs in $\mathcal{G}_{d-1}$. We begin by classifying all $\Gamma\in\mathcal{G}_1$ with $\alpha(\Gamma) < 0$ and all with $\beta(\Gamma)<0$. We then use the results from $\mathcal{G}_1$ to determine all $\Gamma\in\mathcal{G}_2$ with $\alpha(\Gamma) < 0$ and all with $\beta(\Gamma)<0$, continuing similarly for $\mathcal{G}_3$, $\mathcal{G}_4$, and $\mathcal{G}_5$. Finally, we show that no $\Gamma\in \mathcal{G}_d$ has $\alpha(\Gamma) < 0$ or $\beta(\Gamma)<0$ for $d \ge 6$. These results are then summarised in Theorem \ref{recap}, which describes fully all graphs in $\mathcal{G}$ represented by more matrices of rank $n-1$ than rank $n$.

The following outlines some methods used throughout this section.
\begin{itemize}
    \item $\Gamma(t_0,t_1,\dots,t_{d-1},t_d)$ is isomorphic to $\Gamma(t_d,t_{d-1},\dots,t_1,t_0)$ by symmetry. This means that $\alpha(t_0,t_1,\dots,t_{d-1},t_d) = \alpha(t_d,t_{d-1},\dots,t_1,t_0)$, and similar for $\beta$. We use Lemma \ref{alpha-beta} on both forms of $\alpha$ (or $\beta$) to find pairs of conditions which must be satisfied simultaneously, since one form is negative if and only if the other form is negative. These pairs often contradict one another, reducing the number of cases that need to be checked.
    
    \item While $t_0, t_d \ge 0$ for any graph $\Gamma= \Gamma(t_0,t_1,\dots,t_{d-1},t_d)$, all internal values $t_1,\dots,t_{d-1} \ge 1$. This is because the number of edges between the first (or last) vertex in the induced path $P$ of $\Gamma$ and the first (or last) neighbour of the extra vertex $x$ in $P$ may be 0, but the number of edges between any pair of neighbours of $x$ in $P$ must be at least 1. Many of the conditions for $\alpha$ or $\beta$ to be negative require at least one of $t_0$ or $t_d$ to be 0. When determining when $\alpha$ and $\beta$ are negative in $\mathcal{G}_{d}$ using the results from $\mathcal{G}_{d-1}$ and Lemma \ref{alpha-beta}, this often requires an internal value to be 0, which is impossible. This again allows us to reduce the number of cases that need to be checked.
    
    \item If we have shown that $\alpha$ or $\beta$ is equal to an expression involving terms with a power of $-1$ as a coefficient, and we want to prove $\alpha$ or $\beta$ is positive for all values in this case, we may use the fact that it is greater than or equal to the expression resulting from letting all powers of $-1$ be negative simultaneously.

    \item For positive $t_0$ and $t_d$, we note that $\alpha (t_0,t_1,\dots ,t_d)=\alpha (0,t_0,\dots ,t_d,0)$ and $\beta (t_0,t_1,\dots ,t_d)=\beta (0,t_0,\dots ,t_d,0)$. This can be deduced directly from Lemma \ref{alpha-beta} using a symmetry argument, or from the arguments of Section 3. It is a consequence of the fact that the vector $(e_1+e_{n-1})^\top u= 0$ for all $u\in U_{n-1}$, since every $u\in U_{n-1}$ has 1 as its first and last entry. For a graph $\Gamma$ of order $n$ in $\mathcal{G}$, changing the adjacency status of the vertex $x$ with both $x_1$ and $x_{n-1}$ has no effect on the values of $A,B,\alpha$ or $\beta$. In particular, if $\Gamma\in\mathcal{G}_d$ has two vertices $x_1$ and $x_{n-1}$ of degree 1, then the graph $\Gamma'\in\mathcal{G}_{d+2}$ obtained from $\Gamma$ by adding the edges $xx_1$ and $xx_{n-1}$ satisfies $\alpha (\Gamma')=\alpha (\Gamma)$ and $\beta (\Gamma')=\beta (\Gamma)$.

    \item Expressions for $\alpha$ and $\beta$ are found using identities in Lemma \ref{F(p)}. Calculations were done by hand and checked against a SageMath \cite{Sage} program which is included as an appendix.
\end{itemize}

\subsection{Degree 1}
In this subsection, we study graphs of the form $\Gamma(s, t) \in \mathcal{G}_1$ (see Figure \ref{fig Gamma(s,r)}).

\begin{figure}[H]
    \centering
    \includegraphics[width=0.5\linewidth]{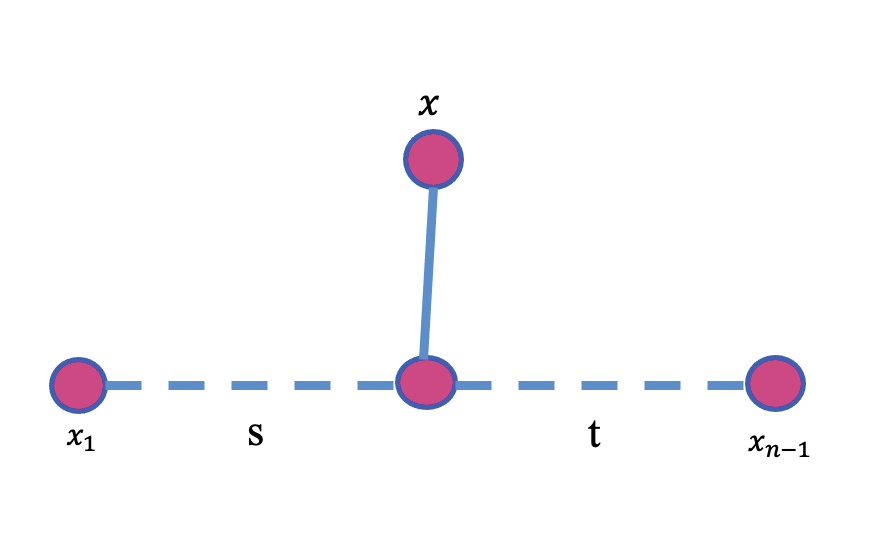}
     \vspace{-1cm}
    \caption{$\Gamma(s,t)$}
    \label{fig Gamma(s,r)}
\end{figure}

Recall $F(n)=\frac{1}{3}(2^{n+1}+(-1)^n)$ and note $F(0)=F(1)=1$, $F(2) = 3$, $F(3) = 5$, $F(4) = 11$, $F(5) = 21$.

\begin{lemma}
\label{base case}
 Let $s$ and $t$ be positive integers.
$$\begin{aligned}
& \alpha(s,t) = 2F(s)F(t) - 2F(s-1)F(t-1)\\
& \beta(s,t)= 8F(s-1)F(t-1) - 2F(s)F(t)\\
\end{aligned}$$
\end{lemma}

\begin{proof}
The number of rank-deficient completions of $M(P_{s+t+1})$ whose right nullspace contains a vector with 1 in position $s+1$ is $\frac{1}{2}A(s,t)$. To form such a completion of $M(P_{s+t+1})$, we may complete the upper left $s\times (s+1)$ region to a matrix $L$ whose rowspace avoids $e_{s+1}^\top$, and complete the lower right $t\times (t+1)$ region to a matrix $R$ whose rowspace avoids the vector $e_1^\top$ in $\F_2^{t+1}$. If $u_L$ and $u_R$ are the unique non-zero vectors in the right nullspaces of $L$ and $R$ respectively, then the last entry of $u_L$ and the first entry of $u_R$ are both 1. There is one way to complete row $s+1$ so that it is orthogonal to the vector $u\in\F_2^{s+t+1}$ whose first $s+1$ and last $t+1$ entries respectively coincide with $u_L$ and $u_R$. By Lemma \ref{last n-1 rows}, the number of choices for $L$ and $R$ are respectively $F(s)$ and $F(t)$, hence $A(s,t)=2F(s)F(t)$.

The number of rank-deficient completions of $M(P_{s+t+1})$ whose right nullspace contains a non-zero vector with 0 in position $s+1$ is $B(s,t)$. To form such a completion, we complete the upper left $s\times s$ region to a completion of $M(P_s)$ of rank $s-1$, complete the lower right $t\times t$ region to a completion of $M(P_t)$ of rank $t-1$, and assign either value to the indeterminate in row $s+1$. The number of choices for the upper left and lower right matrices are respectively $F(s-1)$ and $F(t-1)$, hence $B(s,t)=2F(s-1)F(t-1)$.

 Since $\alpha = A-B$ and $\beta = 4B-A$, this implies the result.
\end{proof}

\begin{theorem}
\label{alpha(s,t)}
$\alpha(s,t)$ is never negative and is zero if and only if $s=t=1.$   
\end{theorem}

\begin{proof}
From Lemma \ref{base case}, $\alpha(s,t)=2 F(s) F(t)-2 F(s-1) F(t-1)$.
Since $F(s)\geq F(s-1)$ and $F(t)\geq F(t-1)$, it follows that $\alpha(\Gamma)$ is always non-negative. Moreover, $\alpha(\Gamma)$ is equal to zero if and only if $s=t=1$, since $F(n) = F(n-1)$ only for $n=1$.
\end{proof}

We now identify when $\beta(s,t)$ is negative, which is used later to determine when $\alpha(r,s,t)<0$.

\begin{theorem}\label{beta(s,t)}
$\beta(s,t)<0$ if and only if $\min (s, t)$ is even.
\end{theorem}
\begin{proof}
From Lemma \ref{base case}, $\beta(s,t)=8 F(s-1) F(t-1)-2 F(s) F(t)$. By Lemma \ref{F(p)}, we derive:
$$\begin{aligned}
\beta(s,t) & =2\left[(-1)^{s+1} F(t) + (-1)^{t+1} F(s)+(-1)^{s+t}\right]
\end{aligned}$$
Therefore $\beta(s,t)$ can be negative, and this happens exactly when $\min(s,t)$ is even. 
\end{proof}

\subsection{Degree 2}
In this subsection, we study graphs of the form $\Gamma(r,s, t) \in \mathcal{G}_2$ (see Figure \ref{fig Gamma(r,s,t)}).

\begin{figure}[h]
    \centering
    \includegraphics[width=0.6\linewidth]{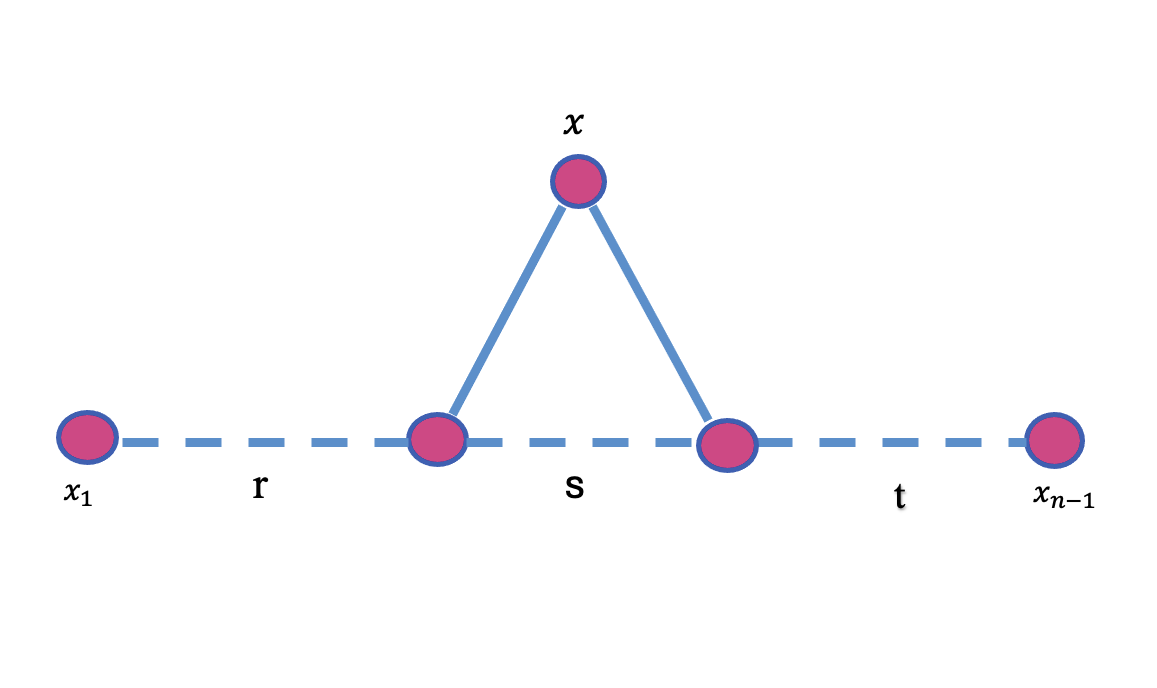}
    \vspace{-1cm}
    \caption{$\Gamma(r,s,t)$}
    \label{fig Gamma(r,s,t)}
\end{figure}

\begin{theorem}
\label{alpha(r,s,t)}
    $\alpha(r,s,t)$ is negative if and only if either:
    \begin{enumerate}[leftmargin=2cm, label = \textbf{Case \arabic*:}]
        \item $r=0$ and $\min(s,t)$ is even, or $t=0$ and $\min(s,r)$ is even
        \item $s=2$ and $\min(r,t)$ is even
    \end{enumerate}
\end{theorem}

\begin{proof}
By Lemma \ref{alpha-beta}, and since $\alpha(r,s,t) = \alpha(t,s,r)$: 
\begin{eqnarray*}
\alpha (r,s,t) & = & \frac{1}{2}F(r)\beta(s,t)+2F(r-1)\alpha (s-1,t) \\
& = & \frac{1}{2}F(t)\beta(s,r)+2F(t-1)\alpha (s-1,r)
\end{eqnarray*}
By Theorem \ref{alpha(s,t)}, the $\alpha$ terms are not negative for any values of $r,s,t$. Therefore $\alpha(r,s,t)$ can only be negative if both $\beta$ terms are negative. By Theorem \ref{beta(s,t)}, this happens if and only if $\min(s,t)$ and $\min(s,r)$ are both even.

\textbf{Case 1:} Suppose $r=0$. Then $\alpha(r,s,t) = \frac{1}{2}F(0)\beta(s,t)+2F(0-1)\alpha (s-1,t) = \frac{1}{2}\beta(s,t)$, which is negative if and only if $\min(s,t)$ is even. Similar for $t=0$.
    
\textbf{Case 2:} Suppose $r,t > 0$. Then $r,s,t \geq 2$ since $\min(s,t)$ and $\min(s,r)$ are both even.

Suppose first that $s\ge 3$. Then $F(s)\ge 5$, $F(r)\ge 3$ and $F(t)\ge 3$. Then $\alpha(r,s,t)$ can be expressed as follows, using the identities in Lemma \ref{F(p)}.

$$\begin{aligned}
\alpha(r,s,t)= \frac{1}{4}\bigg[ 3F(r) F(s) &F(t) +3(-1)^{t+1} F(r) F(s) + 3(-1)^{r+1} F(s) F(t) + 9(-1)^{s+1} F(r)F(t)\\
&+5(-1)^{s+t} F(r)+(-1)^{r+t+1} F(s)+5(-1)^{r+s} F(t)+(-1)^{s+r+t+1}\bigg]\\
\end{aligned}$$

If $s$ is odd, then $r$ and $t$ are both even, giving the following.
$$\begin{aligned}
\alpha(r,s,t) &= \frac{1}{4}\bigg[ 3F(r)F(s)F(t) -3F(r)F(s) -3F(s)F(t) -F(s) + 9F(r)F(t) -5F(r) -5F(t) +1\bigg]\\
&= \frac{1}{4}\bigg[\underbrace{F(r)F(s)F(t)-3F(r)F(s)}_{\ge0\text{ since }F(t)\ge3} + \underbrace{F(r)F(s)F(t)-3F(s)F(t)}_{\ge0\text{ since }F(r)\ge3} +\underbrace{F(r)F(s)F(t)-F(s)}_{\ge 8F(s)\text{ since }F(r),F(t)\ge3}\\
&\phantom{11111111111111111111111} +\underbrace{2F(r)F(t)-5F(r)}_{\ge F(r)\text{ since }F(t)\ge3} +\underbrace{2F(r)F(t)-5F(t)}_{\ge F(t)\text{ since }F(r)\ge3} +5F(r)F(t)+1\bigg]\\
&\ge \frac{1}{4}\big[8F(s) + F(r) + F(t) + 5F(r)F(t) + 1\big] > 0\\
\end{aligned}$$

If $s$ is even and $r$ is odd, then $F(s) \ge 11$, giving the following. \\
$$\begin{aligned}
\alpha(r,s,t) &\ge \frac{1}{4}\bigg[ 3F(r)F(s)F(t) -3F(r)F(s) +3F(s)F(t) - 9F(r)F(t) -5F(r)  -F(s) -5F(t) -1\bigg]\\
&= \frac{1}{4}\bigg[\underbrace{F(r)F(s)F(t)-3F(r)F(s)}_{\ge0\text{ since }F(t)\ge3} + \underbrace{F(r)F(s)F(t)-9F(r)F(t)}_{\ge 2F(r)F(t)\text{ since }F(s)\ge11}+\underbrace{F(s)F(t)-F(s)}_{\ge 2F(s)\text{ since }F(t)\ge3} \\
&\phantom{111111111111111111111111111111} +\underbrace{F(r)F(s)F(t)-5F(r)-5F(t)-1}_{\ge 0\text{ since }F(s)\ge11} +2F(s)F(t)\bigg]\\
&\ge \frac{1}{4}\big[2F(r)F(t) + 2F(s) +2F(s)F(t)\big] > 0\\
\end{aligned}$$

By symmetry, $\alpha(r,s,t)$ is positive if $s$ is even and $t$ is odd. For $s$, $r$, $t$ even, we have the following.
$$\begin{aligned}
\alpha(r,s,t) &= \frac{1}{4}\bigg[ 3F(r)F(s)F(t) -3F(r)F(s) -3F(s)F(t) -9F(r)F(t) +5F(r) -F(s) 
 + 5F(t) -1\bigg]\\
&= \frac{1}{4}\bigg[\underbrace{F(r)F(s)F(t)-3F(r)F(s)}_{\ge0\text{ since }F(t)\ge3} + \underbrace{F(r)F(s)F(t)-3F(s)F(t)}_{\ge0\text{ since }F(r)\ge3} +\underbrace{F(r)-1}_{\ge2\text{ since }F(r)\ge3}\\
&\phantom{1111111111}+\underbrace{\frac{9}{11}F(r)F(s)F(t)-9F(r)F(t)}_{\ge 0\text{ since }F(s)\ge11}+\underbrace{\frac{2}{11}F(r)F(s)F(t)-F(s)}_{\ge \frac{7}{11}F(s)\text{ since }F(r),F(t)\ge3} + 4F(r) + 5F(t)\bigg]\\
&\ge \frac{1}{4}\left[2 + \frac{7}{11}F(s) + 4F(r) + 5F(t)\right] > 0\\
\end{aligned}$$

So $\alpha(r,s,t)>0$ if $r,t>0$ and $s\ge3$. If $s=2$, the formula for $\alpha(r,s,t)$ simplifies as follows.
$$\alpha(r,2,t) = (-1)^{r+1}F(t) + (-1)^{t+1}F(r) + (-1)^{r+t+1}$$
This is negative exactly when $\min(r,t)$ is even. 

\end{proof}


\begin{theorem}
\label{beta(r,s,t)}
$\beta(r, s, t)$ is negative if and only if $s=1$ and $min(r,t)$ are odd. 
\end{theorem}

\begin{proof}
By Lemma \ref{base case} and since $\beta(r,s,t)=\beta(t,s,r)$:
$$
\begin{aligned}
\beta(r,s,t)=&2 F(r) \alpha\left(s,t\right)+2 F(r-1) \beta\left(s-1,t\right)\\
=&2 F(t) \alpha\left(s,r\right)+2 F(t-1) \beta\left(s-1,r\right)
\end{aligned}
$$
By Theorem \ref{alpha(s,t)}, $\alpha(s,t)$ and $\alpha(s,r)$ are never negative. Then $\beta(r,s,t)$ can be negative only if both $\beta(s-1,t)$ and $\beta(s-1,r)$ are negative. By Theorem \ref{beta(s,t)}, this implies $\min(s-1,t)$ and $\min(s-1,r)$ are both even.

If $r=0$, then $
\beta(0, s, t)=2 \alpha(s,t)
$, which is always positive.  Similarly, $\beta(r,s,t) > 0$ if $t=0$.

 If $r,t \ge 1$, we simplify $\beta(r,s,t)$ using $2F(n-1) = F(n) + (-1)^{n+1}$ from Lemma \ref{F(p)} as follows.
$$\begin{aligned} \beta(r, s, t)= & 3 F(r) F(s) F(t)+3(-1)^s F(r) F(t) +(-1)^{r+t} F(s) \\ &-2(-1)^{s+t} F(r)-2(-1)^{s+r} F(t)+(-1)^{s+t+r}\end{aligned}$$
First suppose $s\ge 2$.
Then $s\ge 3$ and $t\ge 2$, since $\min (s-1, t)$ is even. This implies $F(t) \geq 3$ and $F(s)\ge 5$. As a result, $2F(r)F(s)F(t)+3(-1)^{s}F(r)F(t) \geq 7F(r)F(t)$ giving the following.
$$
\beta(r, s, t) \geqslant F(r) F(s) F(t)+7 F(r) F(t)+(-1)^{r+t} F(s)+(-1)^{s+t+r} 
 -2(-1)^{s+t} F(r)-2(-1)^{r+s} F(t)
$$

We note that 
$
7F(r)F(t)>2F(r)+2F(t)+1
$
and $F(r)F(s)F(t)>F(s)$. Hence $\beta(r,s,t)>0$ for $s\ge 2$.

It only remains to consider $s=1$.  We simplify $\beta(r,s,t)$ using Lemma \ref{F(p)} as follows.
$$
\begin{aligned}
\beta(r, 1,t) & =3 F(r) F(t)-3 F(r) F(t)+2(-1)^{t} F(r) +(-1)^{r+t}+2(-1)^{r} F(t)+(-1)^{r+t+1} \\
& =2\left[(-1)^{t} F(r)+(-1)^{r} F(t)\right]
\end{aligned}
$$

Thus $\beta (r,s,t)$ is negative if and only if  $s=1$ and $\min (r,t)$ is odd.

\end{proof}
\newpage
\subsection{Degree 3}
In this subsection, we study graphs of the form $\Gamma(q,r,s, t) \in \mathcal{G}_3$ (see Figure \ref{fig Gamma(q,r,s,t)}).

\begin{figure}[h]
    \centering
    \includegraphics[width=0.75\linewidth]{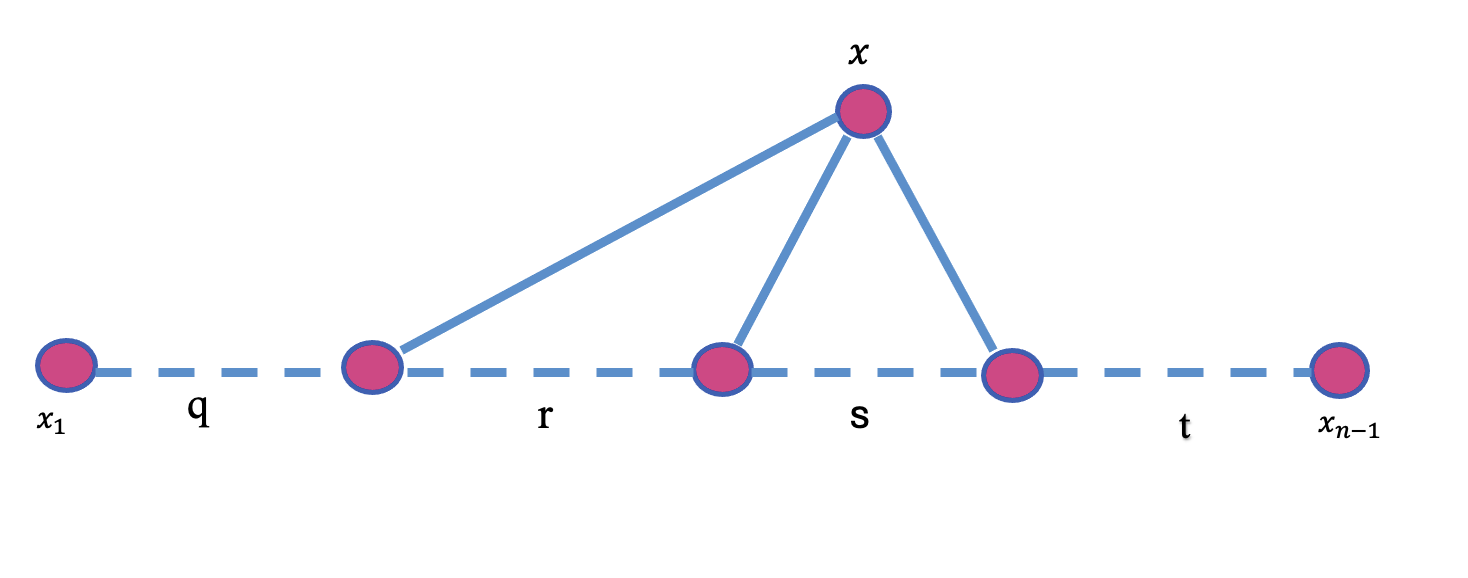}
    \vspace{-1cm}
    \caption{$\Gamma(q,r,s,t)$}
    \label{fig Gamma(q,r,s,t)}
\end{figure}

\begin{theorem}
\label{alpha(q,r,s,t)}
    $\alpha(q,r,s,t)$ is negative if and only if either:
    \begin{enumerate}[leftmargin=2cm, label = \textbf{Case \arabic*:}]
        \item $(q,s)=(0,1)$ and $min (r, t)$ is odd, or $(t,r)=(0,1)$ and $min (q, s)$ is odd.
        \item $s=r=1$ and $min(q, t)$ is even.
    \end{enumerate}
\end{theorem}
\begin{proof}Since $\alpha(q, r, s, t)=\alpha(t, s, r, q), $ we have
\begin{eqnarray*}
    \alpha (q,r,s,t) & = & \frac{1}{2} F(q) \beta(r, s, t)+2 F(q-1) \alpha(r-1, s, t) \\
    & = & \frac{1}{2} F(t) \beta(s, r, q)+2 F(t-1) \alpha(s-1, r, q)
    \end{eqnarray*}

\textbf{Case 1:} If $q=0$, then Lemma \ref{alpha-beta} implies $\alpha(0,r,s,t)=\frac{1}{2} \beta(r,s,t)$. By Theorem \ref{beta(r,s,t)}, this is negative if and only if $s=1$ and $min(r,t)$ is odd. By symmetry, $t=0$ gives the other result in Case 1.

\textbf{Case 2:} If $q,t \geq 1$, Theorems \ref{alpha(r,s,t)} and \ref{beta(r,s,t)} and the above equations for $\alpha(q,r,s,t)$ imply that it can only be negative if both of the following hold.
\begin{itemize}
    \item $\underbrace{s=1 \text{ and } \min(r,t) \text{ is odd}}_{\beta(r,s,t)<0}$, or $\underbrace{r-1=0 \text{ and } \min(s,t) \text{ is even}}_{\alpha(r-1,s,t)<0}$, or $\underbrace{s=2 \text{ and } \min(r-1,t) \text{ is even}}_{\alpha(r-1,s,t)<0}$.
    \item $\underbrace{r=1 \text{ and } \min(s,q) \text{ is odd}}_{\beta(s,r,q)<0}$, or $\underbrace{s-1=0 \text{ and } \min(r,q) \text{ is even}}_{\alpha(s-1,r,q)<0}$, or $\underbrace{r=2 \text{ and } \min(s-1,q) \text{ is even}}_{\alpha(s-1,r,q)<0}$.
\end{itemize}
Therefore we need only consider the cases in which either $r\in\{1,2\}$ or $s\in\{1,2\}$.

\underline{$r=2$:} We simplify $\alpha(q,r,s,t)$ using $2F(n-1) = F(n) + (-1)^{n+1}$ from Lemma \ref{F(p)} as follows.
$$\begin{aligned}
\alpha(q,2,s,t) = \frac{1}{2}\bigg[&12F(q)F(s)F(t) + 3(-1)^{q+1}F(s)F(t)\\
 &+2(-1)^{s+t+1}F(q)+(-1)^{q+t}F(s) + 7(-1)^{q+s}F(t)+3(-1)^{q+s+t+1}\bigg]\\
\end{aligned}$$

If $q=1$, then $\alpha(q,2,s,t)$ simplifies as follows.
$$\begin{aligned}
\alpha(1,2,s,t) &= \frac{1}{2}\bigg[15F(s)F(t) +(-1)^{t+1}F(s) + 7(-1)^{s+1}F(t)+(-1)^{s+t}\bigg]\\
&\ge \frac{1}{2}\bigg[15F(s)F(t) -F(s) -7F(t)-1\bigg]\\
&\ge \frac{1}{2}\bigg[15F(s)F(t) -F(s)F(t) -7F(s)F(t)-F(s)F(t)\bigg]\ge \frac{1}{2}\bigg[6F(s)F(t)\bigg]> 0\\
\end{aligned}$$

If $q\ge2$, then $F(q) \ge 3$, and $\alpha(q,2,s,t)$ simplifies as follows.
$$\begin{aligned}
\alpha(q,2,s,t) = \frac{1}{2}\bigg[&2F(q)F(s)F(t) + 2(-1)^{s+t+1}F(q)\\
 &+10F(q)F(s)F(t)+3(-1)^{q+1}F(s)F(t)+(-1)^{q+t}F(s) + 7(-1)^{q+s}F(t)+3(-1)^{q+s+t+1}\bigg]\\
 \ge \frac{1}{2}\bigg[&\underbrace{2F(q)F(s)F(t) -2F(q)}_{\ge0}+\underbrace{30F(s)F(t)-3F(s)F(t)-F(s)-7F(t)-3}_{\ge 16F(s)F(t)}\bigg] > 0\\
\end{aligned}$$

Therefore $\alpha(q,r,s,t) > 0$ if $r=2$. Similarly, by symmetry, $\alpha(q,r,s,t) > 0$ if $s=2$.

\underline{$r=1$:} We simplify $\alpha(q,r,s,t)$ using $2F(n-1) = F(n) + (-1)^{n+1}$ from Lemma \ref{F(p)} as follows.
$$\begin{aligned}
\alpha(q,1,s,t) = \frac{1}{2}\bigg[&3F(q)F(s)F(t) + 3(-1)^sF(q)F(t) + 3(-1)^{t+1}F(q)F(s)\\
 &+(-1)^{s+t+1}F(q)+2(-1)^{q+t}F(s) + 2(-1)^{q+s}F(t)+2(-1)^{q+s+t+1}\bigg]\\
\end{aligned}$$

We have dealt with $s=2$ previously, so we now separately consider $s=1$ and $s\geq 3$. 

Suppose $s\ge 3$. Since $r=1$, it follows from the bullet points above that for $\alpha (q,r,s,t)<0$, $\min (s,t)$ is even and $\min (s,q)$ is odd. This means that $t\geq 2$. Separating the $3F(q)F(s)F(t)$ term in the expression for $\alpha(q,1,s,t)$ to dominate the potentially negative terms gives the following.

$$\begin{aligned}
2\alpha(q,1,s,t) =& \bigg[F(q)F(s)F(t) + 3(-1)^sF(q)F(t) + 2(-1)^{q+s}F(t) \bigg]\\
&+\bigg[F(q)F(s)F(t) + 3(-1)^{t+1}F(q)F(s)\bigg]+\bigg[\frac{2}{3}F(q)F(s)F(t) + 2(-1)^{q+t}F(s)  \bigg]\\
&+\bigg[\frac{1}{3}F(q)F(s)F(t)+ (-1)^{s+t+1}F(q) + 2(-1)^{q+s+t+1}\bigg]\\
\geq& \bigg[\underbrace{5F(q)F(t) -3F(q)F(t) -2F(t)}_{\geq 0}\bigg]+\bigg[\underbrace{3F(q)F(s) - 3F(q)F(s)}_{\geq 0}\bigg]\\
&+\bigg[\underbrace{2F(q)F(s) - 2F(s)}_{\geq 0} \bigg] +\bigg[\underbrace{5F(q) - F(q) - 2}_{\geq 2F(q)}\bigg]> 0\\
\end{aligned}$$

Similarly, by symmetry, $\alpha(q,r,s,t) > 0$ if $s=1$ and $r \geq 3$.

Now suppose $s=r=1$.
$$\begin{aligned}
\alpha(q,1,1,t) = \frac{1}{2}\bigg[&3F(q)F(t) - 3F(q)F(t) + 3(-1)^{t+1}F(q)\\
 &+(-1)^{t}F(q)+2(-1)^{q+t} + 2(-1)^{q+1}F(t)+2(-1)^{q+t}\bigg]\\
 =(-1&)^{t+1}F(q) + (-1)^{q+1}F(t)+2(-1)^{q+t}\\
\end{aligned}$$
This is negative if and only if $\min(q,t)$ is even.

\end{proof}

\begin{theorem}\label{beta(q,r,s,t)}
    $\beta(q,r,s,t)$ is negative if and only if either:
    \begin{enumerate}[leftmargin=2cm, label = \textbf{Case \arabic*:}]
        \item $(q,s)=(0,2)$ and $min (r, t)$ is even, or $(t,r)=(0,2)$ and $\min (q, s)$ is even.
        \item $q=t=0$, and $\min(r,s)$ is even.
    \end{enumerate}
\end{theorem}

\begin{proof}
Since $\beta(q, r, s, t)=\beta(t, s, r, q), $ we have
\begin{eqnarray*}
    \beta (q,r,s,t) & = & 2 F(q) \alpha(r, s, t)+2 F(q-1) \beta(r-1, s, t) \\
    & = & 2 F(t) \alpha(s, r, q)+2 F(t-1) \beta(s-1, r, q)
    \end{eqnarray*}

If $q=0$, then $\beta(0, r, s, t)=2 \alpha(r, s, t)$ which is negative if $s=2$ and $ min(r,t)$ is even (Case 1), $r=0$ and $\min(s,t)$ is even (not possible since $r \ge 1$), or $t=0$ and $\min(s,r)$ is even (Case 2). Similarly, if $t=0$, then $\beta(q,r,s,0)$ is negative if $r=2$ and $min(s,q)$ is even (finishing Case 1), $s=0$ and $\min(r,q)$ is even (not possible since $s \ge 1$), or $q=0$ and $\min(s,r)$ is even (Case 2).

If $q,t \geq 1$, and since $r,s \geq 1$, Theorems \ref{alpha(r,s,t)} and \ref{beta(r,s,t)} and the above expressions for $\beta(q,r,s,t)$ imply that it can only be negative if both of the following hold.  

\begin{itemize}
    \item $\underbrace{s=2\text{ and }\min(r,t)\text{ is even}}_{\alpha(r,s,t)<0}$, or $\underbrace{s=1\text{ and }\min(r-1,t)\text{ is odd}}_{\beta(r-1,s,t)<0}$.
    \item $\underbrace{r=2\text{ and }\min(s,q)\text{ is even}}_{\alpha(s,r,q)<0}$, or $\underbrace{r=1\text{ and }\min(s-1,q)\text{ is odd}}_{\beta(s-1,r,q)<0}$.
\end{itemize}

We therefore may restrict our attention to $r,s \in \{1,2\}$.

\underline{$r=s=1$:} We simplify $\beta(q,r,s,t)$ as follows.
$$\begin{aligned}
\beta(q,1,1,t) &= 2\big[6F(q)F(t) + (-1)^{q+1}F(t) + (-1)^{t+1}F(q) + (-1)^{q+t}\big]\\
&\ge 2\big[6F(q)F(t) -F(t) -F(q) - 1\big] \ge 2\cdot3F(q)F(t) > 0\\
\end{aligned}$$

\underline{$r=1$, $s=2$:} We simplify $\beta(q,r,s,t)$ as follows.
$$\begin{aligned}
\beta(q,1,2,t) &= 2\big[6F(q)F(t) + 5(-1)^{q+1}F(t) + (-1)^{t}F(q) + (-1)^{q+t+1}\big]\\
&= 2\big[\underbrace{5F(q)F(t) + 5(-1)^{q+1}F(t) + (-1)^{q+t+1}}_{> 0}+\underbrace{F(q)F(t) +  (-1)^{t}F(q)}_{\ge 0}\big]\\
\end{aligned}$$

\underline{$r=s=2$:} We simplify $\beta(q,r,s,t)$ as follows.
$$\beta(q,2,2,t) = 2\big[6F(q)F(t) + 7(-1)^{q+1}F(t) + 7(-1)^{t+1}F(q) + 3(-1)^{q+t}\big]$$

If $q=1$, then $\beta(q,2,2,t)$ simplifies as follows.
$$\beta(1,2,2,t) = 2\big[13F(t) + 10(-1)^{t+1}\big] \ge 2\cdot3F(t) > 0$$

By symmetry, $\beta(q,2,2,1) > 0$.

If $q,t\ge2$, then $F(q), F(t) \ge 3$, and we rearrange $\beta(q,2,2,t)$ as follows.
$$\begin{aligned}
\beta(q,2,2,t) &= 2\big[3F(q)F(t) + 7(-1)^{q+1}F(t) + 3F(q)F(t)+ 7(-1)^{t+1}F(q) +3(-1)^{q+t}\big]\\
&\ge 2\big[9F(t) -7F(t) + 9F(q) -7F(q) -3\big] = 2\big[2F(t) + 2F(q) -3\big] > 0\\
\end{aligned}$$
Therefore $\beta(q,r,s,t) > 0$ for all $q,t \ge 1$.
\end{proof}

\subsection{Degree 4}

In this subsection, we study graphs of the form $\Gamma(p,q,r,s, t) \in \mathcal{G}_4$ (see Figure \ref{fig Gamma(p,q,r,s,t)}).

\begin{figure}[h!]
    \centering
    \includegraphics[width=0.75\linewidth]{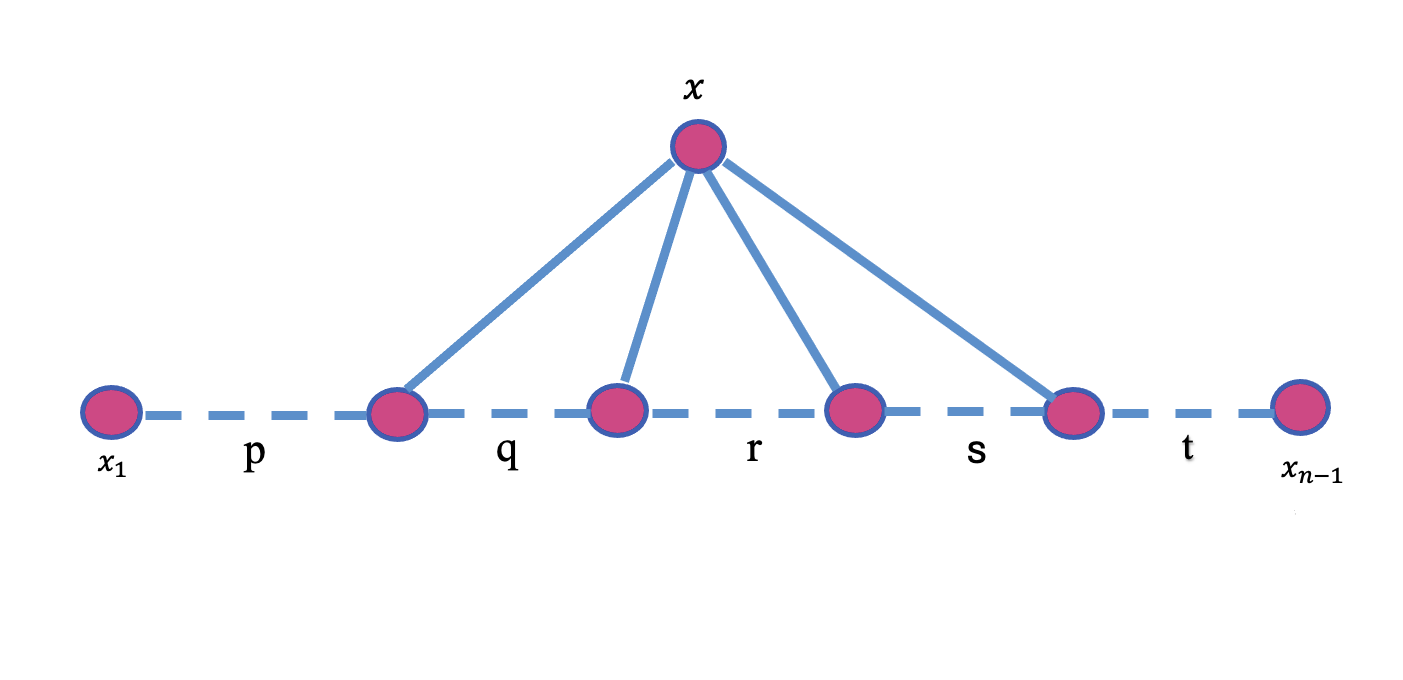}
    \vspace{-1cm}
    \caption{$\Gamma(p,q,r,s,t)$}
    \label{fig Gamma(p,q,r,s,t)}
\end{figure}

\begin{theorem}
\label{alpha(p,q,r,s,t)}
$\alpha(p,q,r,s,t)<0$ if and only if $p=t=0, r=2$ and $\min(q,s)$ is even.
\end{theorem}

\begin{proof}
By Lemma \ref{alpha-beta}, and since $\alpha(p,q,r,s,t)=\alpha(t,s,r,q,p)$:
$$
\begin{aligned}
\alpha(p,q,r, s, t)=&\frac{1}{2} F(p) \beta(q,r,s, t)+2 F(p-1) \alpha(q-1,r,s, t)\\
=&\frac{1}{2} F(t) \beta(s, r, q, p)+2 F(t-1) \alpha(s-1, r, q, p)
\end{aligned}
$$

First suppose $p=0$. Then $\alpha(0,q,r,s,t)=\frac{1}{2} \beta(q,r,s,t)$, so $\alpha(0,q,r,s,t)$ can only be negative if $\beta(q,r,s,t)$ is negative. By Theorem \ref{beta(q,r,s,t)},  and since $q,r,s \ge 1$, this means $p=t=0$, $r=2$, and $\min(q,s)$ is even. Similarly, this is the only case for which $\alpha(p,q,r,s,0)$ is negative. 

Suppose $p,t>0$. Theorems \ref{alpha(q,r,s,t)} and \ref{beta(q,r,s,t)} and the above expressions for $\alpha(p,q,r,s,t)$ imply that it can only be negative if both of the following hold.
\begin{itemize}
    \item $\underbrace{q=s=1 \text{ and } \min(r,t) \text{ is odd}}_{\alpha(q-1,r,s,t)<0}$, or $\underbrace{r=s=1 \text{ and } \min(q-1,t) \text{ is even}}_{\alpha(q-1,r,s,t)<0}$.
   \item $\underbrace{q=s=1 \text{ and } \min(r,p) \text{ is odd}}_{\alpha(s-1,r,q,p)<0}$, or $\underbrace{r=q=1 \text{ and } \min(s-1,p) \text{ is even}}_{\alpha(s-1,r,q,p)<0}$.
\end{itemize}

 From the above cases, $\alpha(p,q,r,s,t)$ can only be negative if $q=s=1$.
$$
\begin{aligned}
   \alpha(p, 1, r, 1, t)=& 3 F(p) F(r) F(t)-3(-1)^r F(p) F(t)-(-1)^{r+t} F(p)-(-1)^{p+r} F(t)-(-1)^{p+t} F(r)
\end{aligned}
$$
If $r=1$, we obtain the following.

$$\begin{aligned}
\alpha(p, 1,1,1, t)&{=6 F(p) F(t)+(-1)^t F(p)+(-1)^p F(t)-(-1)^{p+t}}\\
&\ge 6F(p)F(t) - F(p) - F(t) - 1 \ge 3F(p)F(t) > 0\\
\end{aligned}$$
 If $r > 1$, we simplify $\alpha(p,q,r,s,t)$ as follows.
$$
\begin{aligned}
   \alpha(p, 1, r, 1, t)=& F(p) F(r) F(t)-3(-1)^r F(p) F(t)\\
   &+F(p)F(r)F(t)-(-1)^{r+t} F(p)-(-1)^{p+r}F(t)\\
   &+F(p)F(r)F(t)-(-1)^{p+t} F(r)\\
   \ge& \underbrace{3F(p)F(t)-3F(p)F(t)}_{=0}+\underbrace{3F(p)F(t)-F(p)-F(t)}_{\ge F(p)F(t)}+\underbrace{F(p)F(r)F(t)-F(r)}_{\ge 0}
\end{aligned}
$$
Therefore $\alpha(p,q,r,s,t)$ is negative if and only if $p=t=0$, $r=2$, and $\min(q,s)$ is even.

\end{proof}

\begin{theorem}
\label{beta(p,q,r,s,t)}
    $\beta(p,q,r,s,t)$ is negative if and only if either:
    \begin{enumerate}[leftmargin=2cm, label = \textbf{Case \arabic*:}]
        \item $p=0$, $r=s=1$ and $min (q, t)$ is even, or $t=0$, $q=r=1$ and $\min (p, s)$ is even.
        \item $p=t=0$, $r=1$ and $\min(q,s)$ is odd.
    \end{enumerate}
\end{theorem}
\begin{proof}
Since $\beta(p,q, r, s, t)=\beta(t, s, r, q,p)$, we have
    $$
\begin{aligned}
\beta(p, q, r, s, t)&=2 F(p) \alpha(q, r, s, t)+2 F(p-1) \beta(q-1, r, s, t) \\
& =2 F(t) \alpha(s, r, q, p)+2 F(t-1) \beta(s-1, r, q, p)
\end{aligned}
$$
 If $p=t=0$, since $\alpha(0,q,r,s,0) = \alpha(q,r,s)$, this is negative if and only if $r=1$ and $\min(q,s)$ is odd (Case 2). If $p=0$ and $t>0$, then $\beta(0,q, r, s, t)=2\alpha(q,r,s,t)$. Since $q,r,s\ge1$, this is negative only if $r=s=1$ and $min(q,t)$ is even. Similarly if $t=0$ and $p>0$, then $\beta(p,q,r,s,t)<0$ if and only if $q=r=1$ and $\min(p,s)$ is even (Case 1). 

If $p,t > 0$, then Theorems \ref{alpha(q,r,s,t)} and \ref{beta(q,r,s,t)} and the above expressions for $\beta(p,q,r,s,t)$ imply that it can only be negative if both of the following hold.

\begin{itemize}

\item $\underbrace{r=s=1 \text{ and } \min(q,t) \text{ is even}}_{\alpha(q,r,s,t)<0}$, or $\underbrace{(q,s) = (1,2) \text{ and } \min(r,t) \text{ is even}}_{\beta(q-1,r,s,t)<0}$.
    
\item $\underbrace{r=q=1 \text{ and } \min(s,p) \text{ is even}}_{\alpha(s,r,q,p)<0}$, or $\underbrace{(s,q)=(1,2) \text{ and } \min(r,p) \text{ is even}}_{\beta(s-1,r,q,p)<0}$.

\end{itemize}

 If $r=s=q=1$, then $\min(s,p) = \min(1,p)$ must be even, which contradicts $p,t >0$. If $(r,s,q) \in \{(1,1,2), (1,2,1)\}$, then either $\min(r,t) = \min(1,t)$ or $\min(r,p) = \min(1,p)$ must be even. This again contradicts $p,t >0$, and so $\beta(p,q,r,s,t)$ is not negative for $p,t >0$.
\end{proof}

We conclude that the only graphs $\Gamma \in \mathcal{G}_4$ for which either $\alpha(\Gamma)$ or $\beta(\Gamma)$ is negative arise from the addition of two edges to graphs in $\mathcal{G}_2$, between $x$ and the vertices at the ends of the long induced path.
\newpage
\subsection{Degree 5}

In this subsection, we study graphs of the form $\Gamma(k,p,q,r,s, t) \in \mathcal{G}_5$ (see Figure \ref{fig Gamma(k,p,q,r,s,t)}).
\begin{figure}[h!]
    \centering
    \includegraphics[width=0.75\linewidth]{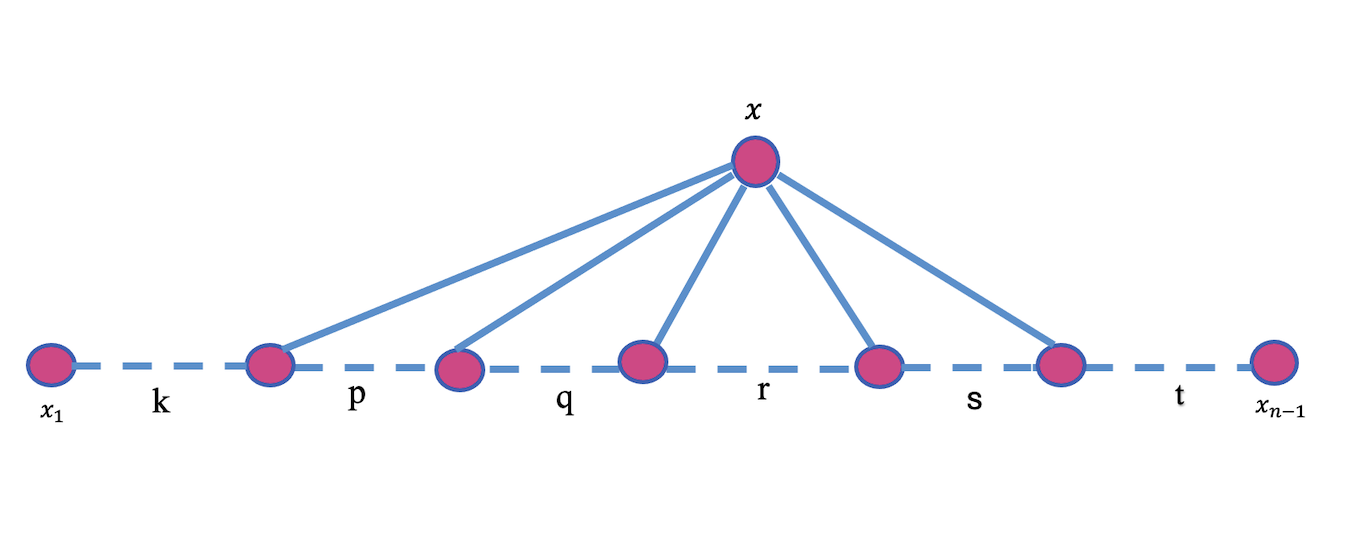}
    \vspace{-0.5cm}
    \caption{$\Gamma(k,p,q,r,s,t)$}
    \label{fig Gamma(k,p,q,r,s,t)}
\end{figure}

\begin{theorem}
\label{alpha(k,p,q,r,s,t)}
    $ \alpha(k, p, q, r, s, t)$ is negative if and only if $k=t=0$, $q=r=1$, and $\min(p,s)$ is even. 
\end{theorem}
\begin{proof}
  
By Lemma \ref{alpha-beta}, and since $\alpha(k, p, q, r, s, t)= \alpha(t, s, r, q, p, k)$:
$$
\begin{aligned}
 \alpha(k, p, q, r, s, t)&=\frac{1}{2} F(k) \beta(p, q, r, s, t)+2 F(k-1) \alpha(p-1, q, r, s, t) \\
&=\frac{1}{2} F(t) \beta(s, r, q, p, k)+2 F(t-1) \alpha(s-1, r, q, p, k)
\end{aligned}
$$

Suppose $k=0$. Then $\alpha(0, p, q, r, s, t)=\frac{1}{2} \beta(p, q, r, s, t)$. By Theorem \ref{beta(p,q,r,s,t)} and since  $p, q, r, s \geq 1$, then $\beta(p, q, r, s, t)$ is negative if and only if $t=0, r=q=1$, and $\min (s, p)$ is even. By symmetry, the same holds if $t=0$.

Next, suppose $k, t > 0$. By Theorems \ref{alpha(p,q,r,s,t)} and \ref{beta(p,q,r,s,t)}, and  since $p,q,r,s \ge 1$, the above expressions for $\alpha(k, p, q, r, s, t)$ imply that it can be negative only if both of the following hold.

\begin{itemize}

\item $\underbrace{t=0, r=q=1 \text{ and } \min(p,s) \text{ is even}}_{\beta(p,q,r,s,t)<0}$, or $\underbrace{ t=0,p=1, r=2 \text{ and } \min(s,q) \text{ is even}}_{\alpha(p-1,q,r,s,t)<0}$.
    
\item $\underbrace{k=0, r=q=1 \text{ and } \min(s,p) \text{ is even}}_{\beta(s,r,q,p,k)<0}$, or $\underbrace{ k=0, s=1, q=2\text{ and } \min(r,p) \text{ is even}}_{\alpha(s-1,r,q,p,k)<0}$.

\end{itemize}

 $\beta(p,q,r,s,t)<0$ implies $q=1$ and $\alpha(s-1,r,q,p,k)<0$ implies $q=2$, and therefore they cannot hold simultaneously. Similarly, $\beta(s,r,q,p,k)<0$ and $\alpha(p-1,q,r,s,t)$ cannot hold simultaneously due to contradicting requirements for $r$. If both $\alpha(p-1, q, r, s, t)$ and $\alpha(s-1, r, q, p)$ are negative, $\min (s, q)=\min (1,2)$  is not even, another contradiction.

If $\beta(p,q,r,s,t)<0$ and $\beta(s,r,q,p,k)<0$, this means that $k=t=0$. Therefore, $$
\alpha(0, p, q, r, s, 0)=\frac{1}{2} \beta(p, q, r, s, 0)=\frac{1}{2} \beta(0, s, r, q, p)=\alpha(s, r, q, p)=\alpha(p, q, r, s),
$$
which is negative if and only if $q=r=1$ and $\min (p, s)$ is even.
\end{proof}

\begin{theorem}
\label{beta(k, p, q, r, s, t)}
 $\beta(k, p, q, r, s, t)$  is positive for all $k,p,q,r,s,t \in \mathbb{N}$.
\end{theorem}
\begin{proof}
 By Lemma \ref{alpha-beta}, and since $\beta(k, p, q, r, s, t)= \beta(t, s, r, q, p, k)$: 
    $$
\begin{aligned}
 \beta(k, p, q, r, s, t)&=2 F(k) \alpha(p, q, r, s, t)+2 F(k-1) \beta(p-1, q, r, s, t) \\
&=2 F(t) \alpha(s, r, q, p, k)+2 F(t-1 ) \beta(s-1, r, q, p, k)
\end{aligned}
$$

If $k=0$, then $\beta(0, p, q, r, s, t)=2 \alpha(p, q, r, s, t)$. By Theorem \ref{alpha(p,q,r,s,t)}, and since  $p\ge1$,  it follows that $\alpha(p, q, r, s, t)$ is never negative. By symmetry,  $\beta(k, p, q, r, s, t)$ is never negative if $t=0$.

Suppose $k,t>0$. By Theorems \ref{alpha(p,q,r,s,t)} and \ref{beta(p,q,r,s,t)}, and since $p,q,r,s \ge 1$, the above expressions for $\beta(k,p,q,r,s,t)$ imply that it can be negative only if both of the following hold.

\begin{itemize}
\item $\underbrace{p=r=s=1 \text{ and } \min(q,t) \text{ is even}}_{\beta(p-1,q,r,s,t)<0}$.
\item $\underbrace{p=q=s=1 \text{ and } \min(r,k) \text{ is even}}_{\beta(s-1,r,q,p,k)<0}$.
\end{itemize}

But since $r=q=1$ in this case, $\min(q,t)=1$ and $\min(r,k)=1$, which contradicts the requirement that both are even. Therefore $\beta(k,p,q,r,s,t)>0$ for all $k,p,q,r,s,t \in \mathbb{N}$.
\end{proof}

We conclude that the only graphs $\Gamma \in \mathcal{G}_5$ for which $\alpha(\Gamma)$ is negative arise from the addition of two edges to graphs in $\mathcal{G}_3$, between $x$ and the vertices at the ends of the long induced path.

\subsection{Degree $\geq 6$}
In this final subsection, we prove that any graph of the form $\Gamma\left(t_0, t_1, \ldots, t_d\right) \in \mathcal{G}_d$ is represented by more matrices of rank $n$ than rank $n-1$, for $d\ge6$.

\begin{theorem}
    $\alpha(t_0,t_1,\cdots, t_d)$ and $\beta(t_0,t_1,\cdots, t_d)$ are positive if $d \geq 6$. 
\end{theorem}
\begin{proof}
Induction on $d$. From Lemma \ref{alpha-beta}, we have the following.
$$
\begin{aligned}
 \alpha(t_0,t_1,\dots,t_d)&=\frac{1}{2} F(t_0) \beta(t_1,t_2,,\dots,t_d)+2 F(t_0-1) \alpha(t_1-1,t_2,\dots,t_d) \\
\beta(t_0,t_1,\dots,t_d)&=2 F(t_0) \alpha(t_1,t_2,\dots,t_d)+2 F(t_0-1) \beta(t_1-1,t_2,\dots,t_d) \\
\end{aligned}
$$

Suppose $d=6$. Theorem \ref{beta(k, p, q, r, s, t)} implies that $\beta(t_1,t_2,\cdots, t_d)$ and $\beta(t_1-1,t_2,\cdots, t_d)$ are never negative. Since $t_1 \ge 1$, Theorem \ref{alpha(k,p,q,r,s,t)} implies that $\alpha(t_1-1,t_2\cdots, t_d)$ and $\alpha(t_1,t_2\cdots, t_d)$ are never negative. Therefore $\alpha(t_0,t_1,\cdots, t_d)$ and $\beta(t_0,t_1,\cdots, t_d)$ are positive for $d=6$.

Now suppose $\alpha(t_1,t_2,\cdots, t_d)$ and $\beta(t_1,t_2,\cdots, t_d)$ are positive for some $d\ge6$. By the induction hypothesis, each term on the right-hand side of the above equations is positive. Therefore $\alpha(t_0,t_1,\cdots, t_d)$ and $\beta(t_0,t_1,\cdots, t_d)$ are positive if $d \geq 6$.
\end{proof}
\newpage
\section{Conclusions}
We recap our main results in the following theorem.

\begin{theorem}
\label{recap}
The following are all graphs in $\mathcal{G}$ represented by more matrices of rank $n-1$ than rank $n$.

\begin{itemize}
    \item $\Gamma(0,s,t)$ with $\min(s,t)$ even.
    \item $\Gamma(s,2,t)$ with $\min(s,t)$ even.
    \item $\Gamma(0,s,1,t)$ with $\min(s,t)$ odd.
    \item $\Gamma(s,1,1,t)$ with $\min(s,t)$ even.
    \item $\Gamma(0,s,2,t,0)$ with $\min(s,t)$ even.
    \item $\Gamma(0, s, 1, 1, t, 0)$ with $\min(s,t)$ even. 
\end{itemize}
\end{theorem}

\begin{figure}[h!]

    \centering
    \includegraphics[width=0.94\linewidth]{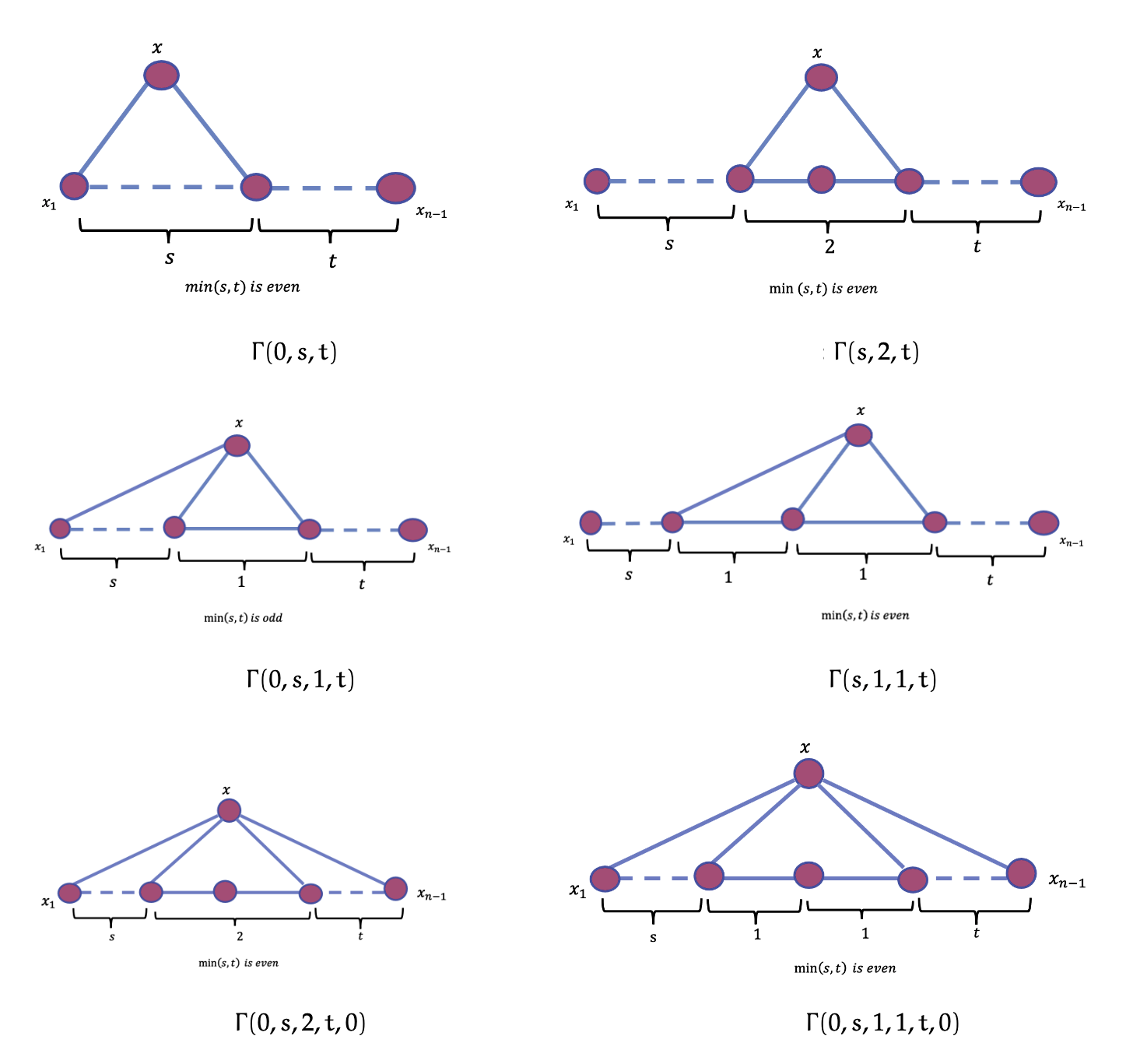}
    \caption{Graphs represented by more matrices of rank $n-1$}
    \label{fig Gammas}
\end{figure}

While $\mathcal{G}$ contains many examples of graphs $\Gamma$ on $n$ vertices for which $R_{n-1}(\Gamma) > R_n(\Gamma)$, many such graphs exist outside of $\mathcal{G}$. Identifying classes of graphs with this property remains an interesting area of research. We have identified all such graphs containing a long induced cycle, and this will be the subject of a forthcoming paper. Beyond these classes of graphs, the question remains open, and it would be of particular interest to exhaustively identify all $\Gamma$ with $R_{n-1}(\Gamma) > R_n(\Gamma)$.

\section*{Acknowledgements}
The research of the first author was supported by a Ph.D. scholarship from Saudi Government Scholarship, in cooperation with Tabuk University.



{
\section*{Appendix}
\vspace{-1.5cm}\hspace{-2cm}\includegraphics[width=1.3\textwidth]{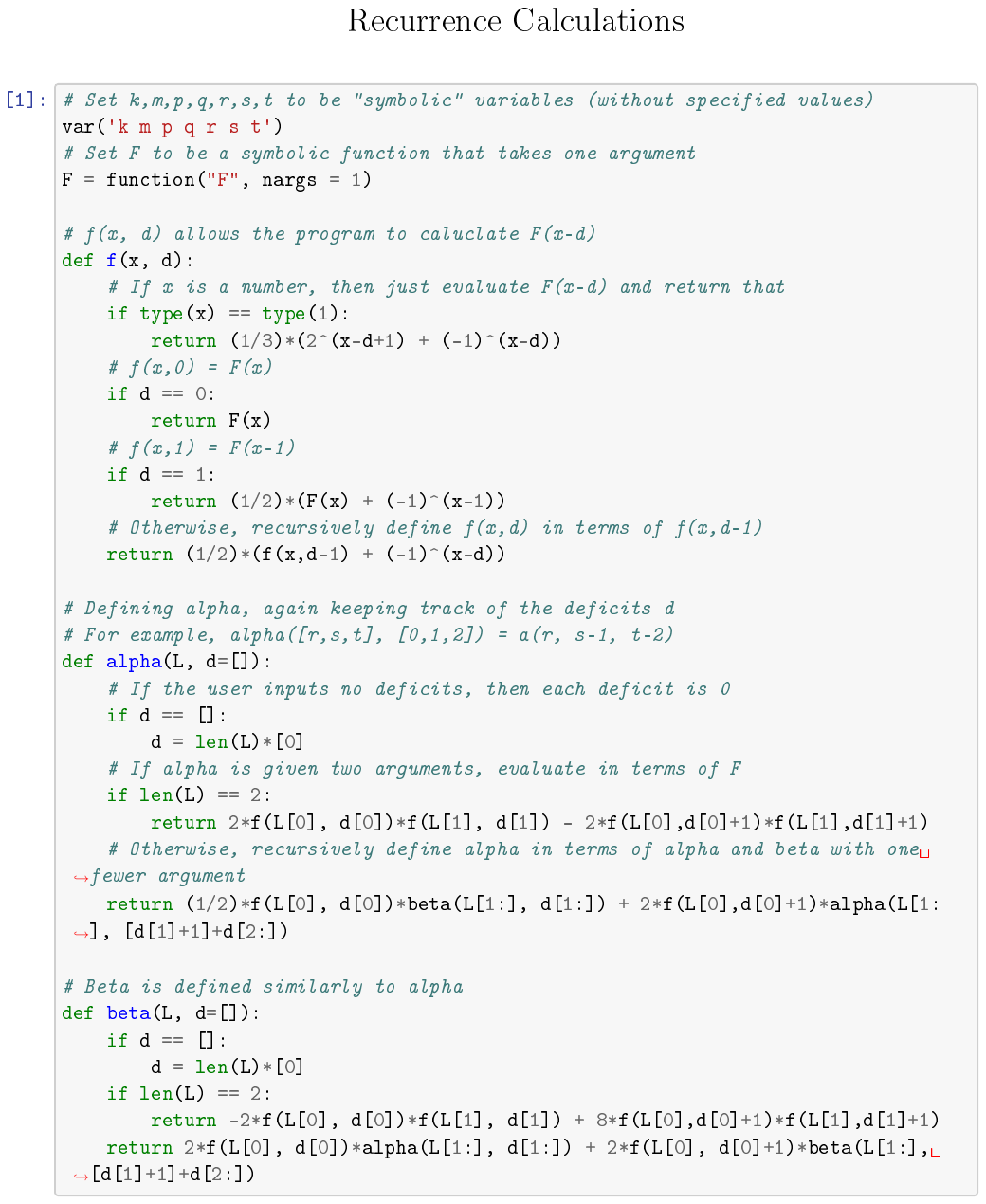}
\vspace{-2cm}\hspace{-2cm}\includegraphics[width=1.3\textwidth]{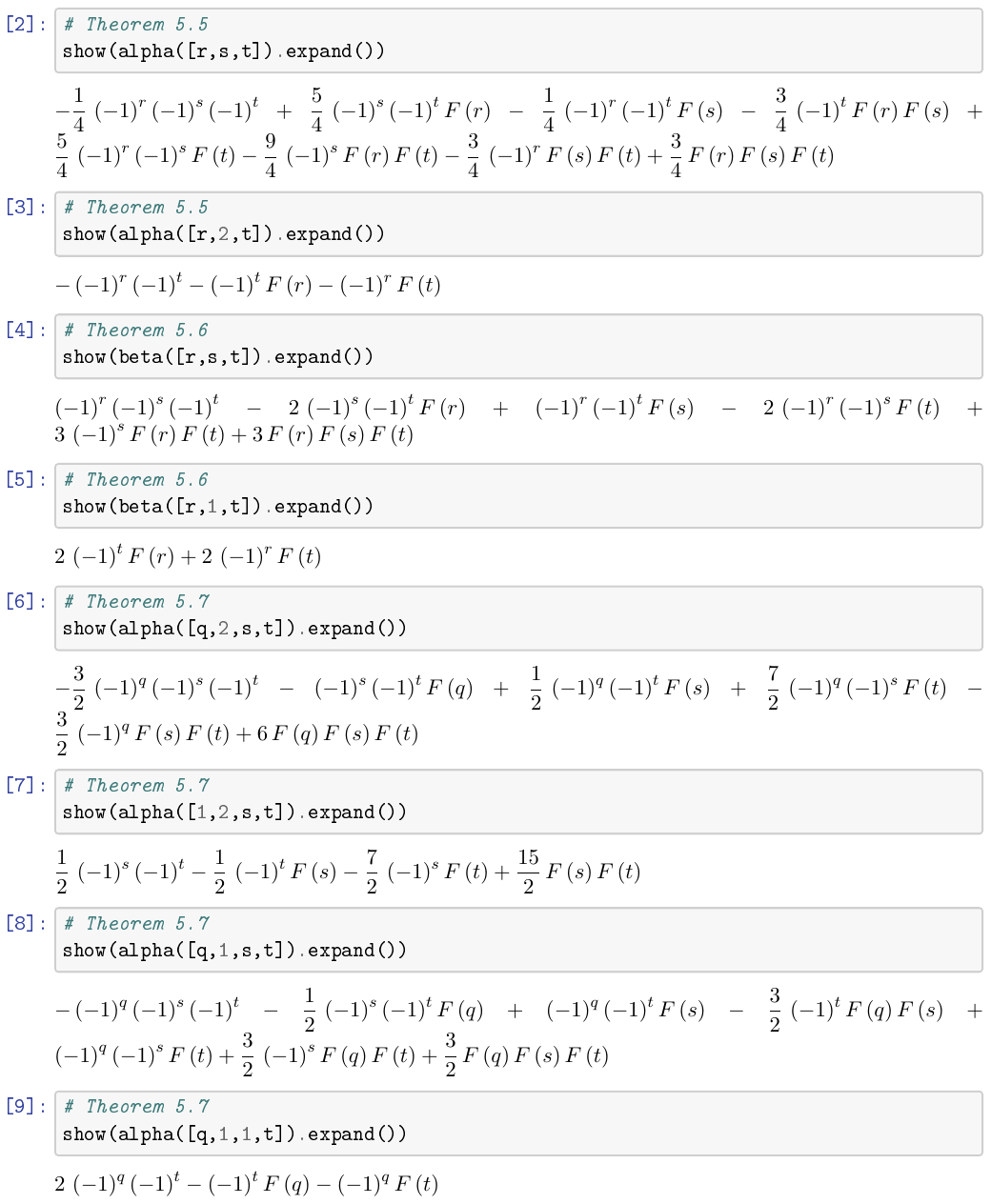}
\vspace{-2cm}\hspace{-2cm}\includegraphics[trim = {0 5cm 0 0}, width=1.3\textwidth]{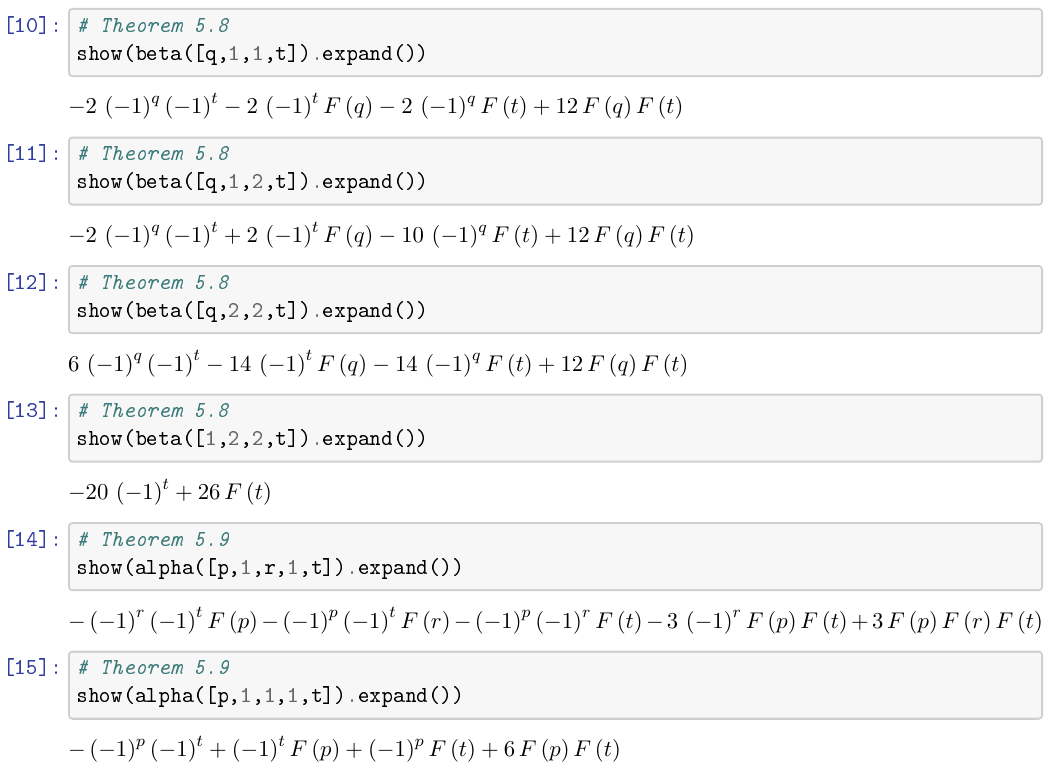}
}
\end{document}